\documentclass[12pt]{iopart}

\usepackage{iopams}  
\usepackage{lipsum}
\usepackage{amsfonts}
\usepackage{graphicx}
\usepackage{epstopdf}
\usepackage{algorithmic}

\usepackage{stackrel}
\usepackage[shortlabels]{enumitem}
\usepackage{epsfig}
\usepackage[final]{pdfpages}
\usepackage{float}
\usepackage{amssymb} 
\usepackage{enumitem}
\usepackage{comment}

\newcommand{\fft}{\mathcal{F}}
\newcommand{\paren}[1]{\left(#1\right)}
\newcommand{\parenn}[1]{\left[#1\right]}
\newcommand{\parennn}[1]{\left\{#1\right\}}
\newcommand{\parennnn}[1]{\left|#1\right|}
\newcommand{\norm}[1]{\left\|#1\right\|}

\newcommand{\Hbb}{\mathcal{H}}
\newcommand{\Nbb}{\mathbb{N}}
\newcommand{\Rbb}{\mathbb{R}}

\newcommand{\set}[2]{\left\{#1\,\left|\,#2\right.\right\}}

\newcommand{\mmap}[3]{#1:\,#2\rightrightarrows #3\,}

\usepackage{amsopn}

\newcommand*{\ee}{\mathtt{e}}
\newcommand*{\ii}{\mathtt{j}}
\newcommand*{\jj}{\mathtt{j}}

\newcommand{\FFT}{\operatorname{\mathcal{F}}}

\newcommand{\Abs}[1]{\left\vert#1\right\vert}
\newcommand*{\F}{\operatorname{\mathcal{F}}}

\newcommand{\Hcal}{{\mathcal H}}

\newcommand{\Id}{\mathtt{Id}}
\newcommand{\dist}{\mathtt{dist}}
\newcommand{\Fix}{\mathtt{Fix\,}}

\newtheorem{thm}{Theorem}[section]

\newtheorem{fact}[thm]{Fact}
\newtheorem{remark}[thm]{Remark}
\newtheorem{lemma}[thm]{Lemma}
\newtheorem{definition}[thm]{Definition}
\newtheorem{example}[thm]{Example}
\newtheorem{proposition}[thm]{Proposition}
\newtheorem{corollary}[thm]{Corollary}
\newtheorem{theorem}[thm]{Theorem}

\usepackage[textsize=small
]{todonotes}

\begin{document}

\title[Projection methods for high-NA phase retrieval]{Projection methods for high numerical aperture phase retrieval}

\author{Nguyen Hieu Thao${}^{1,2,*}$, Oleg Soloviev${}^{1,3}$, Russell Luke${}^4$ and Michel Verhaegen${}^1$}

\address{${}^1$ Delft Center for Systems and Control, Delft University of Technology, 2628CD Delft, The Netherlands.\\
${}^2$ Department of Mathematics, School of Education, Can Tho University, Can Tho, Vietnam.\\
${}^3$ Flexible Optical B.V., Polakweg 10-11, 2288 GG Rijswijk, The Netherlands.\\
${}^4$ Institut f\"ur Numerische und Angewandte Mathematik, Universit\"at G\"ottingen, 37083 G\"ottingen, Germany.
}
\ead{h.t.nguyen-3@tudelft.nl}
\vspace{10pt}

\begin{abstract}
We develop for the first time a mathematical framework in which the class of projection algorithms can be applied to high numerical aperture (NA) phase retrieval.
Within this framework, we first analyze the basic steps of solving the high-NA phase retrieval problem by projection algorithms and establish the closed forms of all the relevant prox-operators.
We then study the geometry of the high-NA phase retrieval problem and the obtained results are subsequently used to establish convergence criteria of projection algorithms in the presence of noise.
Making use of the vectorial point-spread-function (PSF) is, on the one hand, the key difference between this paper and the literature of phase retrieval mathematics which deals with the scalar PSF.
The results of this paper, on the other hand, can be viewed as extensions of those concerning projection methods for low-NA phase retrieval.
Importantly, the improved performance of projection methods over the other classes of phase retrieval algorithms in the low-NA setting now also becomes applicable to the high-NA case.
This is demonstrated by the accompanying numerical results which show that available solution approaches for high-NA phase retrieval are outperformed by projection methods.
\end{abstract}

\vspace{2pc}


 

\section{Introduction}\label{s:intro}

\emph{Phase retrieval} is an important inverse problem in optics which aims at recovering a complex signal at the pupil plane of an optical system given a number of intensity measurements of its Fourier transform.
It appears in many scientific and engineering fields, including microscopy \cite{Arr99,KimZho16}, astronomy imaging \cite{DaiFie87,HarTho00}, X-ray crystallography \cite{Har93,Mil90}, adaptive optics \cite{AntVer15,VisBruVer16,VisVer13,MugBlaIdi06}, \emph{etc.}
For optical systems with low numerical aperture (NA), a vast number of phase retrieval algorithms have been devised, for example, in \cite{BauComLuk02,CanEldStrVor13,DoeThaVer18,Fie82,Fie13,GerSax72,Hau86,Luk05,LukSabTeb19,SheEldCohChaMiaSeg15,NguLukSolVer20,NguSolVer21} based on the \emph{Fresnel approximation} stating that the intensity distribution in the focal plane and the complex signal in the pupil plane are related via the Fourier transform \cite{Goo05}.
Among solution approaches for low-NA phase retrieval, the widely used class of projection algorithms, which can be viewed as descendants of the classical Gerchberg-Saxton algorithm \cite{GerSax72}, outperforms the other classes by almost every important performance measure: computational complexity, convergence speed, accuracy and robustness \cite[page 410]{LukSabTeb19}.
For high-NA optical systems, the vector nature of light cannot be neglected and point-spread-functions (PSFs) are formed according to a more involved imaging formulation \cite{LinRodKenLarDai12,Man09,McC02,RicWol59}, which is called the \emph{vectorial PSF} to be distinguished from the \emph{scalar} one according to the Fresnel diffraction equation.
In contrast to low-NA settings, only few solution algorithms have been proposed for phase retrieval in high-NA settings \cite{BraDirJanHavNes05,HanGusAgaSed03,NguSolVer20}. 

In this paper, we develop for the first time a mathematical framework in which the class of projection algorithms can be applied to high-NA phase retrieval.
Within this framework, we first analyze the basic steps of solving the high-NA phase retrieval problem by projection algorithms and establish the closed forms of all the relevant prox-operators.
We then study the geometry of the high-NA phase retrieval problem and the obtained results are subsequently used to establish convergence criteria of projection algorithms in the presence of noise.
Making use of the vectorial PSF is, on the one hand, the key difference between this paper and the literature of phase retrieval mathematics which mostly deals with the scalar PSF, see, for example,  \cite{BauComLuk02,Fie82,GerSax72,Gon82,Luk08,LukBurLyo02,SouThiSchFerCouUns16,Sou77,NguLukSolVer20,NguSolVer21}.
The results of this paper, on the other hand, can be viewed as extensions of those concerning projection methods for low-NA phase retrieval.
Importantly, the improved performance of projection methods over the other classes of phase retrieval algorithms in the low-NA setting \cite[page 410]{LukSabTeb19} now also becomes applicable to the high-NA case.
This is demonstrated by the accompanying numerical results which show that all available solution approaches for high-NA phase retrieval are outperformed by projection methods.

The remainder of this paper is organized as follows.
Mathematical notation is introduced in section \ref{subs:math} and the vectorial PSF model is recalled in section \ref{subs:vector PSF}.
In section \ref{s:problem formulation}, several feasibility models of the high-NA phase retrieval problem are formulated based on the vectorial PSF model (data fidelity) and the prior knowledge of the solutions.
In section \ref{s:prox-operator}, closed forms of the projectors on the constituent sets of the feasibility models are established.
In section \ref{s:projection algorithm}, we discuss projection algorithms for solving the feasibility problems in section \ref{s:problem formulation}.
Section \ref{s:prox-regular} is devoted to studying the geometry of the high-NA phase retrieval problem where the constituent sets of feasibility are proven to be \emph{prox-regular} at the points relevant for the subsequent convergence analysis.
Section \ref{s:convergence analysis} is devoted to analyzing convergence of projection algorithms for solving the high-NA phase retrieval in the presence of noise.
As the first ingredient of convergence, the \emph{pointwise almost averagedness} property of projection algorithms is established in section \ref{subs:a.a.} based on the prox-regularity of the component sets proven in section \ref{s:prox-regular}.
The second condition of convergence concerning the mutual arrangement of the component sets around the solution \cite{Kru06,KruLukTha18} is beyond the analysis of this paper.
Convergence criteria are formulated in section \ref{subs:convergence results}.
Numerical simulations are presented in section \ref{s:numerical result}.

\subsection{Mathematical notation}\label{subs:math} 

The underlying space in this paper is a finite dimensional Hilbert space 
denoted by $\mathcal{H}$.
The Frobenius norm denoted by $\|\cdot\|$ is used for both vector and array objects.
Equality, inequalities and mathematical operations such as the multiplication, the division, the square, the square root, the amplitude $|\,\cdot\,|$, the argument $\arg(\cdot)$ and the real part $\Re(\cdot)$ acting on arrays are understood element-wise.
The imaginary unit is $\jj = \sqrt{-1}$.
The distance function associated to a set $\Omega\subset \mathcal{H}$
is defined by
\begin{equation*}
	\dist(\cdot, \Omega) \colon \mathcal{H} \to\Rbb_+\colon x\mapsto \inf_{w\in\Omega}\norm{x-w},
\end{equation*}
and the set-valued mapping
\begin{equation}\label{def:projector}
	\mmap{P_\Omega}{\mathcal{H}}{\Omega}\colon
	x\mapsto \set{w\in \Omega}{\norm{x-w}=\dist(x,\Omega)}
\end{equation}
is the corresponding \emph{projector}.
A selection $w\in P_\Omega(x)$ is called a {\em projection} of $x$ on $\Omega$.
When the projection $w$ is unique, we write $P_\Omega(x)=w$ instead of $P_\Omega(x)=\{w\}$ for brevity.
The \emph{reflector} associated with $\Omega$ is accordingly defined by $R_\Omega \equiv 2P_\Omega-\Id$, where $\Id$ is the \emph{identity mapping}.
Since only projections on either affine or compact sets are involved in the analysis of this paper, the existence of projections is guaranteed.
The fixed point set of a self set-valued mapping $T:\mathcal{H}\rightrightarrows \mathcal{H}$ is defined by $\Fix T \equiv \left\{x\in \mathcal{H} \mid x\in T(x)\right\}$, see, for example, \cite[Definition 2.1]{LukNguTam18}.
An iterative sequence $x_{k+1}\in T(x_k)$ generated by $T$ is said to \emph{converge linearly} to a point $x$ with rate $c\in (0,1)$ if there exists a number $\gamma>0$ such that
\[
\norm{x_k-x} \le \gamma c^k\quad \forall k\in \Nbb.
\]
For $x\in \mathcal{H}$, $\Omega\subset \mathcal{H}$ and an integer $m\ge2$, we make use of the following notation
\begin{equation}\label{duplication}
	[x]_m \equiv \underbrace{(x,x,\ldots,x)}_{m~\mbox{times}}\; \mbox{ and }\;
	[\Omega]_m \equiv \left\{[w]_m \mid w\in \Omega\right\}.
\end{equation}
Our other basic notation is standard; cf. \cite{DonRoc14,Mor06,VA}.
The open ball with radius $\delta>0$ and center $x$ is denoted by $\mathbb{B}_\delta(x)$.

\subsection{Vectorial point-spread-functions}\label{subs:vector PSF}

This section presents the imaging formulation considered in the paper.
For high-NA optical systems, PSFs should be modeled according to the \emph{vector diffraction theory}, see, for example, \cite{HanGusAgaSed04,LinRodKenLarDai12,Man09,McC02,RicWol59}.
More specifically, the $x,y,z$ components of the electromagnetic field right after the lens should be considered separately for the $x$ and $y$ components of the electromagnetic field just before the lens.
Here we consider collimated beams and hence the $z$ component of the field before the lens is zero.
Let the unit electromagnetic fields in the $x$ and $y$ directions just before the lens respectively produce the fields right after the lens with components denoted by $\paren{E_\mathtt{XX}(x,y), E_\mathtt{XY}(x,y), E_\mathtt{XZ}(x,y)}$ and $\paren{E_\mathtt{YX}(x,y), E_\mathtt{YY}(x,y), E_\mathtt{YZ}(x,y)}$, where 
$(x,y)$ are the coordinates in the lens aperture denoted by $\mathcal{P}$.
Let the lens aperture $\mathcal{P}$ be normalized to have radius equal the ${\rm NA}$ value.
Then according to, for example, \cite[Table 3.1]{Man09}, the latter functions are given by
\begin{eqnarray}
	E_\mathtt{XX}(x,y) = \, 1-\frac{k_\mathtt{X}^2(x,y)}{1+k_\mathtt{Z}(x,y)},\;
	&E_\mathtt{YX}(x,y) = \, -\frac{k_\mathtt{Y}(x,y) k_\mathtt{X}(x,y)}{1+k_\mathtt{Z}(x,y)},\nonumber\\
	E_\mathtt{XY}(x,y) = -\frac{k_\mathtt{X}(x,y) k_\mathtt{Y}(x,y)}{1+k_\mathtt{Z}(x,y)},\;
	&E_\mathtt{YY}(x,y) = 1 -\frac{k_\mathtt{Y}^2(x,y)}{1+k_\mathtt{Z}(x,y)},\label{polarisation}
	\\
	E_\mathtt{XZ}(x,y) = - k_\mathtt{X}(x,y),
	&E_\mathtt{YZ}(x,y) = - k_\mathtt{Y}(x,y),\nonumber
\end{eqnarray}	
where $(k_\mathtt{X}(x,y),k_\mathtt{Y}(x,y),k_\mathtt{Z}(x,y))$ is the \emph{unit wave vector} determined for each point $(x,y)$ of the lens aperture $\mathcal{P}$ and satisfies
\begin{equation*}
\max_{(x,y)\in \mathcal{P}} \paren{k_\mathtt{X}^2(x,y) + k_\mathtt{Y}^2(x,y)} = \max_{(x,y)\in \mathcal{P}} \paren{x^2+y^2} = {\rm NA}^2,
\end{equation*}
where the maximum is attainable on the boundary of $\mathcal{P}$ and $\textrm{NA}$ is the NA value.
In particular, the following equality will be used frequently in our subsequent analysis:
\begin{equation}\label{sum=2}
	\sum_{c\,\in \mathcal{I}}E_{c}^2(x,y) = 2,\quad \forall (x,y)\in \mathcal{P},
\end{equation}
where and elsewhere in the paper, the letter $c$ stands for elements of the index set:
\begin{equation}\label{index I}
	\mathcal{I}\equiv \left\{\mathtt{XX}, \mathtt{XY}, \mathtt{XZ}, \mathtt{YX}, \mathtt{YY}, \mathtt{YZ}\right\}.
\end{equation}
In the sequel, the coordinates $(x,y)$ of two-dimensional arrays objects will be dropped for brevity, for example, we simply write $E_{c}$ instead of $E_{c}(x,y)$.

Each of the right-hand side terms in (\ref{polarisation}) can be treated as a corresponding amplitude modulation in the entrance pupil for calculation of a PSF according to the Fresnel diffraction equation:
\begin{equation}\label{psf}
	p_{c}(\mathcal{A},\Phi) = \left|\F\left(E_{c} \cdot \mathcal{A} \cdot \ee^{\ii \Phi}\right)\right|^2,\quad (\forall c\in \mathcal{I})
\end{equation}
where $\mathcal{A}$ and $\Phi$ are respectively the amplitude and phase of the collimated beam in the pupil plane, and $\F$ is the two-dimensional Fourier transform.
The six constituent PSFs according to (\ref{psf}) then can be used to calculate the vectorial PSF corresponding to any linear polarization of light in the entrance pupil.
For unspecified polarization state of light, they are added incoherently as follows:
\begin{equation*}
	I(\mathcal{A},\Phi) = \sum_{c\,\in \mathcal{I}} p_{c}(\mathcal{A},\Phi).
\end{equation*}
Thus, the vectorial PSF with an additional \emph{phase diversity} $\phi_d$ is accordingly given by
\begin{equation}\label{vector PSF}
	I(\mathcal{A},\Phi,\phi_d) = \sum_{c\,\in \mathcal{I}} \Abs{\F \paren{E_{c} \cdot \mathcal{A} \cdot \ee^{\ii\paren{\Phi+\phi_d}}}}^2.
\end{equation}

\begin{figure}[b!]
\centering
\includegraphics[width=0.3\linewidth]{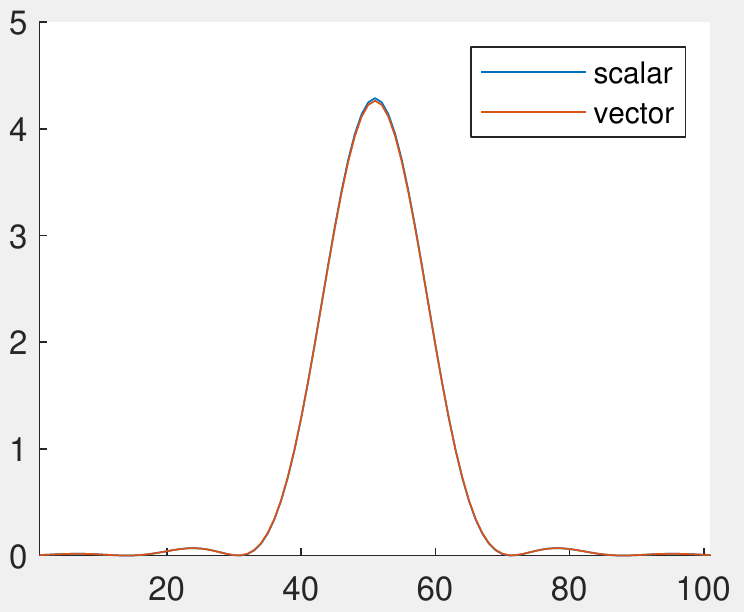}\quad
\includegraphics[width=0.31\linewidth]{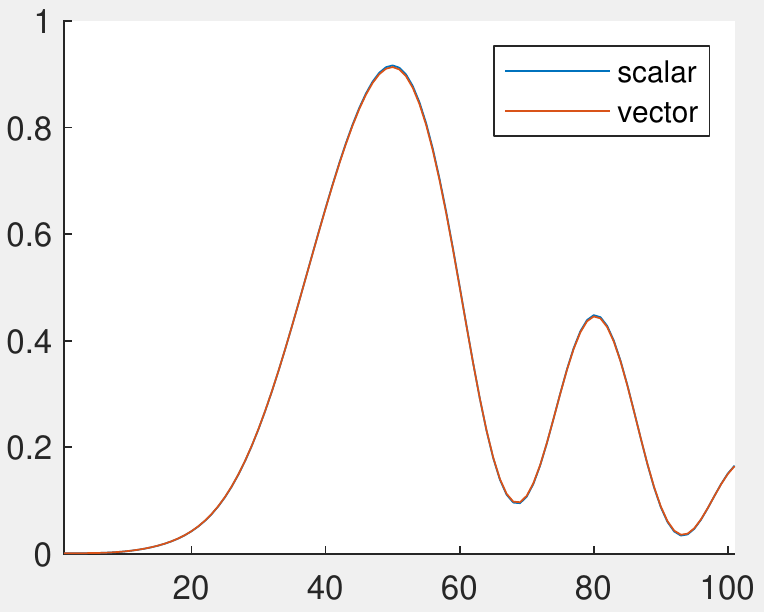}\vspace*{.5cm}\\
\includegraphics[width=0.3\linewidth]{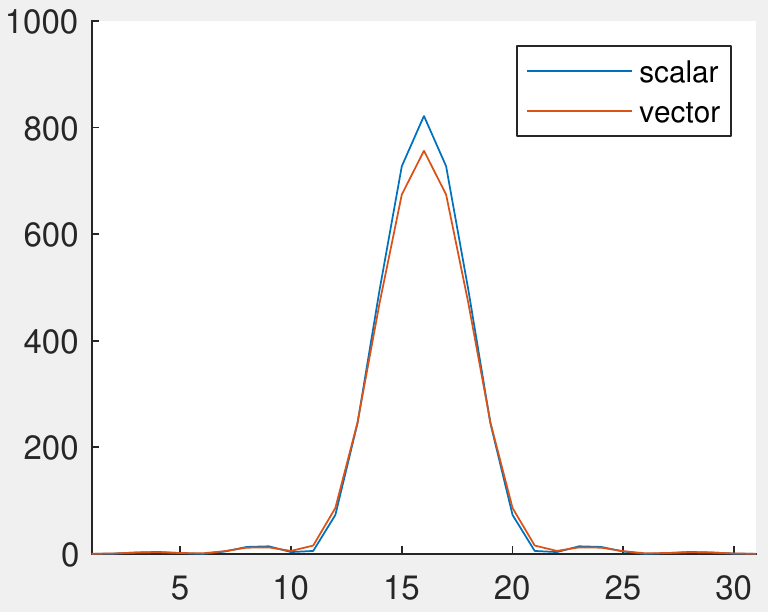}\quad
\includegraphics[width=0.3\linewidth]{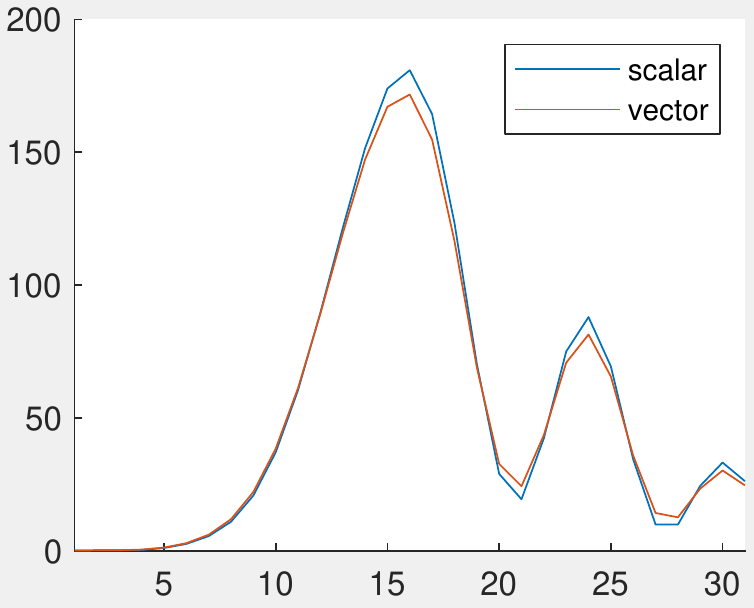}\vspace*{.5cm}\\
\includegraphics[width=0.3\linewidth]{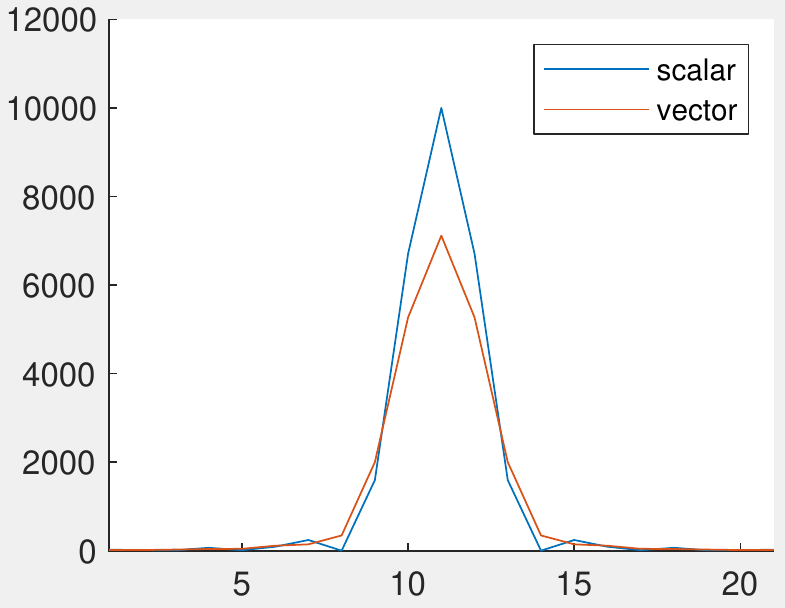}\quad
\includegraphics[width=0.3\linewidth]{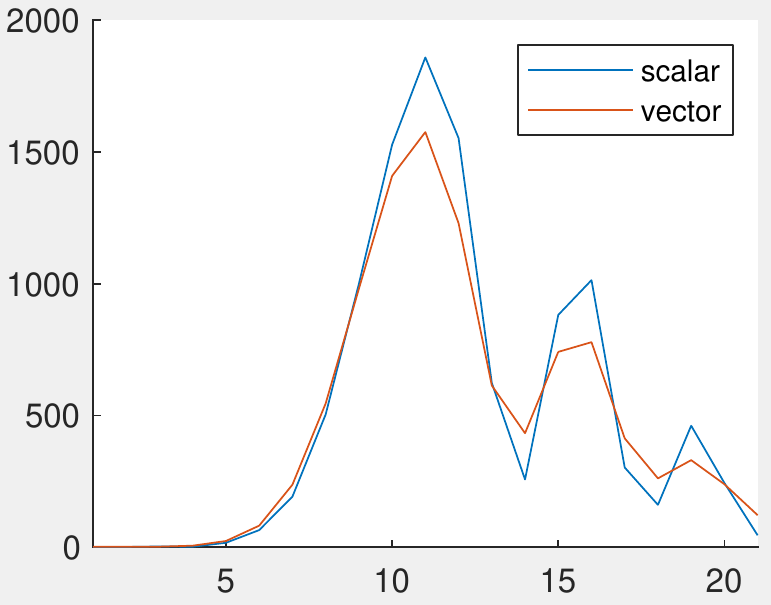}
\caption{Comparison between the \textit{scalar} and \textit{vectorial} PSF models for various NA values -- 0.15, 0.55 and 0.95 in order from top to bottom. In each plot, a pair of corresponding cross-sections of the scalar (in blue color) 
and vectorial (in red color) PSFs are shown. The PSFs on the left-hand-side are without phase aberration ($\Phi=0$) and the ones on the right-hand-side are with phase aberration ($\Phi\neq 0$). The two PSF models differ more for higher NA values and the discrepancy becomes substantial for NA from 0.55.}
\label{fig:S_vs_V_comparison}
\end{figure}

The computational complexity of the vectorial PSF model (\ref{vector PSF}) as a sum of six constituent components is approximately six times higher than the one of the scalar PSF.
There is hence a trade-off between computational complexity and model accuracy in choosing the imaging model of high-NA phase retrieval.
Let us briefly analyze this matter.
Figure \ref{fig:S_vs_V_comparison} reports a short comparison between the \textit{scalar} and \textit{vectorial} PSF models for various NA values -- 0.15, 0.55 and 0.95 in order from top to bottom.
The left-hand-side column of Figure \ref{fig:S_vs_V_comparison} shows PSFs without phase aberration and the second one shows those with phase aberration.
In each plot, a pair of corresponding cross-sections of the scalar (the blue curves) and vectorial (the red curves) PSFs are shown.
It is clear that for low-NA values ($0.15$), the use of the vectorial PSF is superfluous as the two models are almost identical while the vectorial one is much more computationally expensive.
The scalar and vectorial PSF models differ more for higher-NA imaging systems and for particular application purposes their discrepancy can become substantial for NA values from 0.55.

\section{Problem formulation}\label{s:problem formulation}

\subsection{High-NA phase retrieval}\label{subs:PR}

This paper considers the same setting of high-NA phase retrieval as in \cite{NguSolVer20}.
For an unknown phase aberration $\Phi \in \mathbb{R}^{n\times n}$, let $r_d \in \mathbb{R}_+^{n\times n}$ be the measurement of $m$ PSF images $I(\mathcal{A},\Phi,\phi_d)$ generated by (\ref{vector PSF}) with phase diversities $\phi_d$ $\paren{d=1,2,\ldots,m}$.
The \emph{high-NA phase retrieval} problem is to restore $\Phi$ given $r_d$ and $\phi_d$ $\paren{d=1,2,\ldots,m}$ as well as the physical parameters of the optical system.
Mathematically, we consider the problem of finding $\Phi \in \mathbb{R}^{n\times n}$ such that
\begin{equation}\label{phase problem}
	r_d = \sum_{c\,\in \mathcal{I}}\left|\FFT\paren{E_{c}\cdot \mathcal{A} \cdot  \ee^{\jj\paren{\Phi+\phi_d}}}\right|^2  + w_d \quad \paren{d=1,2,\ldots,m},
\end{equation}
where $\mathcal{A}$ is the possibly unknown amplitude of the generalized pupil function (GPF) and $w_d\in \mathbb{R}^{n\times n}$ $\paren{d=1,2,\ldots,m}$ represent the discrepancies between the theoretically predicted data and the actually measured one, for example, due to noise and model deviations.

\subsection{Feasibility models}\label{subs:FP}

In this section, we formulate feasibility models of the phase retrieval problem (\ref{phase problem}) in two scenarios of application -- known and unknown amplitude of the GPF.
According to the vectorial PSF (\ref{vector PSF}), we consider the underlying space
\[
\Hbb \equiv \underbrace{\mathbb{C}^{n\times n}\times \mathbb{C}^{n\times n}\times \dots \times \mathbb{C}^{n\times n}}_{6 \mbox{~times}}.
\]
In the sequel, for each $\paren{x_\mathtt{XX}, x_\mathtt{XY}, x_\mathtt{XZ}, x_\mathtt{YX}, x_\mathtt{YY}, x_\mathtt{YZ}} \in \Hbb$ and $z \in \mathbb{C}^{n\times n}$, we make use of the following notation in accordance with (\ref{index I}):
\begin{eqnarray*}
	\paren{x_{c}}_{c\,\in\mathcal{I}} &\equiv \paren{x_\mathtt{XX}, x_\mathtt{XY}, x_\mathtt{XZ}, x_\mathtt{YX}, x_\mathtt{YY}, x_\mathtt{YZ}},
	\\
	\paren{x_{c} \cdot z}_{c\,\in\mathcal{I}} &\equiv \paren{x_\mathtt{XX} \cdot z, x_\mathtt{XY} \cdot z, x_\mathtt{XZ} \cdot z, x_\mathtt{YX} \cdot z, x_\mathtt{YY} \cdot z, x_\mathtt{YZ} \cdot z}.
\end{eqnarray*}

\subsubsection{Unknown GPF amplitude.}
The following set captures the first constraint of a solution to (\ref{phase problem}) as an element of $\Hbb$:
\begin{equation}\label{Omega_0}
	\Omega_0 \equiv  \parennn{\paren{E_{c}\cdot x}_{c\,\in\mathcal{I}} \in \Hbb \;\mid\; x\in \mathbb{C}^{n\times n}},
\end{equation}
where the six matrices $E_{c}$ ($c\,\in\mathcal{I}$) are defined in (\ref{polarisation}).
Note that $\Omega_0$ is linear subspace of $\Hbb$ with $\dim(\Omega_0)$ being one sixth of $\dim(\Hbb)$.
For $d=1,2,\ldots,m$, the intensity constraint set corresponding to phase diversity $\phi_d$ is given by
\begin{equation}\label{Omega_d}
\Omega_d \equiv  \parennn{\paren{x_{c}}_{c\,\in\mathcal{I}} \in \Hbb \;\mid\; \sum_{c\,\in\mathcal{I}}\left|\FFT\paren{E_{c}\cdot x_{c}\cdot \ee^{\jj\phi_d}}\right|^2 = r_d}.
\end{equation}
The high-NA phase retrieval problem (\ref{phase problem}) then can be addressed via the following $(m+1)$-set feasibility:
\begin{equation}\label{PR1}
	\mbox{find}\quad x\in \bigcap_{d=0}^m \Omega_d.
\end{equation}
The following two-set feasibility models formulated in the product spaces, which are equivalent to (\ref{PR1}) in the case of consistent feasibility (\textit{i.e.}, the intersection is nonempty) \cite{Pie84}, are widely used in practice:
\begin{eqnarray}\label{PR3}
	&\mbox{find }\quad u \in A\cap B \subset \Hbb^m,
	\\\label{PR4}
	&\mbox{find }\quad u \in D \cap B^+ \subset \Hbb^{m+1},
\end{eqnarray}
where
\begin{eqnarray}
	A \equiv  \parennn{(x,x,\ldots,x) \in \Hbb^m \mid x\in \Omega_0},\;
	B \equiv \Omega_1\times \Omega_2\times \cdots\times \Omega_m,\label{A&B}
	\\
	D \equiv  \parennn{(x,x,\ldots,x) \in \Hbb^{m+1} \mid x\in \Hbb},\;
	B^+ \equiv  \Omega_0\times \Omega_1\times \cdots\times \Omega_m.\nonumber
\end{eqnarray}

\subsubsection{Known GPF amplitude.} When the amplitude $\mathcal{A}$ of the GPF is known, it brings stronger constraint on the solutions of (\ref{phase problem}) than (\ref{Omega_0}):
\begin{equation}\label{chi}
	\chi \equiv  \left\{\paren{E_{c}\cdot \mathcal{A} \cdot \ee^{\jj\Phi}}_{c\,\in \mathcal{I}} \in \Hbb \mid \Phi \in \mathbb{R}^{n\times n}\right\}.
\end{equation}
Similar to the case of unknown GPF amplitude, the phase retrieval problem (\ref{phase problem}) then can be addressed via one of the following feasibility models:
\begin{eqnarray}
	\label{PR1'}
	\mbox{find}\quad &x\in \chi \cap \Omega_1 \cap \Omega_2 \cap \cdots \cap \Omega_m,
	\\
	\label{PR3'}
	\mbox{find}\quad &u \in A_{\chi}\cap B,
	\\
	\label{PR4'}
	\mbox{find}\quad &u \in D\cap B_\chi,
\end{eqnarray}
where
\begin{equation}\label{A_chi & B_chi}
	A_{\chi} \equiv  \parennn{(x,x,\ldots,x) \in \Hbb^m \mid x\in \chi},\;
	B_\chi \equiv  \chi \times \Omega_1\times \Omega_2\times \cdots\times \Omega_m.
\end{equation}

\begin{remark}[inconsistent feasibility]\label{r:inconsistent}
Due to noise and model deviations, the intersections in (\ref{PR1}), (\ref{PR3}), (\ref{PR4}), (\ref{PR1'}), (\ref{PR3'}) and (\ref{PR4'}) are empty for all practical purposes.
Keeping in mind, however, that a projection algorithm as a fixed point operator built on a feasibility model is not limited to finding points of intersection but convergence of its iterations to a fixed point of the associated operator is desirable and sufficient in all scenarios of feasibility.
Such fixed points should admit interpretation in terms of meaningful (approximate) solutions to the practical problem  captured by the feasibility model.
We refer the reader to \cite[page 414]{LukSabTeb19} and \cite[Remark 5]{NguSolVer21} for more details on inconsistency of feasibility formulations of (low-NA) phase retrieval.
\end{remark}

\begin{remark}[effectiveness of the feasibility approach]\label{r:feasibility for PR}
It was observed in the recent benchmark paper \cite[page 410]{LukSabTeb19} concerning low-NA phase retrieval that algorithms built on feasibility models outperform all other classes of solution methods by almost every important performance measure.
This observation has strongly encouraged the current work which extends this class of algorithms for high-NA phase retrieval.
\end{remark}

\begin{remark}[choice of feasibility models]\label{r:choose FP}
Depending on specific setting of phase retrieval, one feasibility model can result in better approximate solutions than another.
\end{remark}

\section{Calculation of projectors}\label{s:prox-operator}

The decisive step of solving feasibility problems is to calculate the projectors associated to the relevant sets.
The results of this section, which can be viewed as the high-NA extensions of the ones concerning prox-operators for low-NA phase retrieval \cite{BauComLuk02,Fie82,Gon82,Luk08,LukBurLyo02,SouThiSchFerCouUns16}, enable us to address the feasibility models formulated in section \ref{subs:FP} using projection algorithms.

For convenience let us first introduce further notation and preliminary results.
For each $d=1,2,\ldots,m$ we define the operator $M_d:\Hbb \to \Hbb$ by
\begin{equation}\label{M_d}
	x = \paren{x_{c}}_{c\,\in \mathcal{I}}\;\; \mapsto\;\; M_d(x) \equiv \paren{\FFT\paren{x_{c}\cdot \ee^{\jj\phi_d}}}_{c\,\in \mathcal{I}},
\end{equation}
which is a unitary transform and its inverse is given by
\begin{equation}\label{M_d^-1}
	M_d^{-1}: x = \paren{x_{c}}_{c\,\in \mathcal{I}}\;\; \mapsto\;\; \paren{\fft^{-1}\paren{x_{c}}\cdot\ee^{-\jj\phi_d}}_{c\,\in \mathcal{I}}.
\end{equation}

We then define the matrix-valued function $\mathcal{G}_d: \mathcal{H} \to \mathbb{R}_+^{n\times n}$ by
\begin{equation}\label{G_d}
	x = \paren{x_{c}}_{c\,\in \mathcal{I}} \in \mathcal{H} \mapsto \mathcal{G}_d(x) \equiv \sqrt{\sum_{c\,\in \mathcal{I}}\left|M_d(x)_{c}\right|^2} \in \mathbb{R}_+^{n\times n}.
\end{equation}

\begin{fact}[continuity of $\mathcal{G}_d$]\label{f:continuous of G_d}
The matrix-valued function $\mathcal{G}_d$ is continuous on $\Hcal$.
\end{fact}
\textbf{Proof.}
Since compositions of continuous mappings are continuous, the statement follows from the continuity of $M_d$ and the elementwise amplitude and summation operations.
\hfill$\square$
\smallskip

We define the set $S_d \subset \Hbb$ by
\begin{equation} \label{S_d}
	S_d \equiv \parennn{x = \paren{x_{c}}_{c\,\in \mathcal{I}} \in \Hbb \;\mid\; \sum_{c\,\in \mathcal{I}}\left|x_{c}\right|^2 = r_d}.
\end{equation}

In the sequel, we also make use of the following set of indices:
\[
\mathcal{J} \equiv  \{\xi = (\xi_1,\xi_2) \mid 1\le \xi_1,\xi_2 \le n\},
\]
and for any $x = \paren{x_{c}}_{c\,\in \mathcal{I}} \in \Hbb$, we denote $x[\xi] \equiv \paren{x_{c}[\xi]}_{c\,\in \mathcal{I}}$ and $S_d[\xi] \equiv \{x[\xi] \mid x\in S_d\}$.
In other words, the index of discretized two-dimensional signals (for example, $x_c$ for each $c\,\in \mathcal{I}$) is specified by $\xi$ in square brackets while the index of higher-dimensional arrays such as $x = \paren{x_{c}}_{c\,\in \mathcal{I}} \in \Hbb$ or $S_d$ is defined inductively.

\begin{fact}[projection on $S_d$]\label{f:P_S_d}
It holds that
\begin{equation}\label{P_S_d=prod}
	P_{S_d}(z) = \prod_{\xi \in \mathcal{J}} P_{S_d[\xi]}\paren{z[\xi]}\quad (\forall\, z\in \Hbb),
\end{equation}
where
\begin{equation}\label{P_S_d(xi)}
	P_{S_d[\xi]}\paren{z[\xi]} = \cases{\frac{\sqrt{r_d[\xi]}}{\norm{z[\xi]}} \cdot z[\xi] & if $\|z[\xi]\| \neq 0$,\\
	\parennn{s \in \mathbb{C}^6 \;\mid\; \norm{s}^2 = r_d[\xi]} & if $\|z[\xi]\| = 0$.}
\end{equation}
\end{fact}
\textbf{Proof.}
The product structure (\ref{P_S_d=prod}) is inherent from the product structure of the set $S_d$, that is,
\[
    S_d = \prod_{\xi \in \mathcal{J}} S_d[\xi]\; \mbox{ and hence }\;
    P_{S_d}(z) = \prod_{\xi \in \mathcal{J}} P_{S_d[\xi]}\paren{z[\xi]}\quad (\forall\, z\in \Hbb).
\]
Let us compute $P_{S_d[\xi]}\paren{z[\xi]}$ for each index $\xi \in \mathcal{J}$.
Note that by its definition the set $S_d[\xi]$ is the sphere in $\mathbb{C}^6$ centered at the origin with radius $\sqrt{r_d[\xi]}$, that is,
\begin{equation}\label{S_d_xi}
	S_d[\xi] = \parennn{s \in \mathbb{C}^6 \;\mid\; \norm{s}^2 = r_d[\xi]}.
\end{equation}
Hence its associated projector $P_{S_d[\xi]}$ admits the closed form (\ref{P_S_d(xi)}) as claimed.
\hfill$\square$

The next two results are widely known in the literature of feasibility analysis \cite{Pie84}.
Recall the notation in (\ref{duplication}).

\begin{fact}[projection on diagonals, $P_D$]\label{f:P_D}
For any $w = \paren{w_0,w_1,\ldots,w_m}\in \Hbb^{m+1}$ it holds that
\[
P_D(w) = \parenn{\overline{w}}_{m+1}\; \mbox{ with }\; \overline{w} \equiv \frac{1}{m+1}\sum_{d=0}^m w_d.
\]
\end{fact}

\begin{fact}[projection on product sets, $P_B$]\label{f:P_B}
For any $w = \paren{w_1,w_2,\ldots,w_m}\in \Hbb^m$ it holds that
\begin{equation*}
	P_B(w) = \prod_{d=1}^m P_{\Omega_d}(w_d).
\end{equation*}
\end{fact}

We can now calculate the projectors associated with the sets defined in section \ref{subs:FP}.

\begin{lemma}[projection on $\Omega_0$]\label{l:P_Omega_0}
For any $x = \paren{x_{c}}_{c\,\in \mathcal{I}} \in \Hbb$ it holds that
\begin{equation*}
	P_{\Omega_0}(x) = \paren{E_{c}\cdot z}_{c\,\in\mathcal{I}} \;\mbox{ with }\; z \equiv \frac{1}{2}\sum_{c\,\in \mathcal{I}}\paren{E_{c}\cdot x_{c}}.
\end{equation*}
\end{lemma}
\textbf{Proof.}
By definition (\ref{Omega_0}) of $\Omega_0$ and the definition of projector in (\ref{def:projector}), $\paren{E_{c}\cdot a\cdot \ee^{\jj\Psi}}_{c\,\in \mathcal{I}}$ is a projection of $x$ on $\Omega_0$ if and only if $\paren{a,\Psi}$ is a solution to the following minimization problem:
\begin{equation}\label{min problem_Omega_0}
	\min_{a \in \mathbb{R}_+^{n\times n}, \Psi \in \mathbb{R}^{n\times n}} \quad  \norm{x-\paren{E_{c}\cdot a\cdot \ee^{\jj\Psi}}_{c\,\in \mathcal{I}}}^2.
\end{equation}
The objective function of (\ref{min problem_Omega_0}) can be rewritten as
\begin{eqnarray*}
\norm{x-\paren{E_{c}\cdot a\cdot \ee^{\jj\Psi}}_{c\,\in \mathcal{I}}}^2
&= \|x\|^2 + \|\paren{E_{c}\cdot a}_{c\,\in \mathcal{I}}\|^2 
\\
&- 2\Re\paren{\paren{\sum_{c\,\in \mathcal{I}}\paren{E_{c}\cdot a\cdot x_{c}}} \cdot \ee^{-\jj\Psi}}.    
\end{eqnarray*}
The problem (\ref{min problem_Omega_0}) is hence equivalent to the following one:
\begin{equation}\label{max problem_Omega_0}
	\min_{a \in \mathbb{R}_+^{n\times n}, \Psi \in \mathbb{R}^{n\times n}} \quad \|\paren{E_{c}\cdot a}_{c\,\in \mathcal{I}}\|^2 - 2\Re\paren{\paren{\sum_{c\,\in \mathcal{I}} \paren{E_{c}\cdot a\cdot x_{c}}} \cdot \ee^{-\jj\Psi}}.
\end{equation}
The structure of (\ref{max problem_Omega_0}) allows us to solve for $\Psi$ and $a$ successively though its objective function is not completely separable in $a$ and $\Psi$.
Indeed, since $a$ has no influence on the argument of $\sum_{c\,\in \mathcal{I}}\paren{E_{c}\cdot a\cdot x_{c}}$, the set of optimal $\Psi$ is given by
\begin{equation}\label{Psi opt}
	\left\{\Psi \in \mathbb{R}^{n\times n} \mid \Psi \in \arg\paren{\sum_{c\,\in \mathcal{I}} \paren{E_{c}\cdot x_{c}}}\right\}.
\end{equation}
Plugging the optimal $\Psi$ above into (\ref{max problem_Omega_0}), we arrive at minimizing a quadratic function of variable $a$.
Taking into account that $\sum_{c\,\in \mathcal{I}} \parennnn{E_{c}}^2 = 2J_n$ by (\ref{sum=2}) where $J_n$ is the all-ones matrix of size $n\times n$, we obtain by direct calculation that the unique optimal $a$ is given by
\[
a = \frac{\parennnn{\sum_{c\,\in \mathcal{I}} \paren{E_{c}\cdot x_{c}}}}{\sum_{c\,\in \mathcal{I}} \parennnn{E_{c}}^2} = \frac{1}{2}\parennnn{\sum_{c\,\in \mathcal{I}} \paren{E_{c}\cdot x_{c}}}.
\]
Note that for any index $\xi\in \mathcal{J}$, if $a[\xi]=0$ , then $\Psi[\xi]$ does not play any role in the product $a[\xi]\ee^{\jj\Psi[\xi]}$.
Otherwise, $\Psi[\xi]$ is uniquely determined in view of (\ref{Psi opt}).
Hence, the unique optimal solution to (\ref{min problem_Omega_0}) is given by
\[
z = a\cdot \ee^{\jj\Psi} = \frac{1}{2}\sum_{c\,\in \mathcal{I}}\paren{E_{c}\cdot x_{c}}.
\]
The proof is complete.
\hfill$\square$

\begin{lemma}[projection on $\Omega_d$]\label{l:P_Omega_d}
For each $d=1,2,\ldots,m$ and any $x \in \Hbb$ it holds that
\begin{equation*}
	P_{\Omega_d}(x) = M_d^{-1}(y)
\end{equation*}
where $y \in \Hbb$ is characterized as follows.
\begin{enumerate}[(i)]
\item\label{l:P_Omega_d i} If $\mathcal{G}_d(x)[\xi] \neq 0$, then $y[\xi] = \frac{\sqrt{r_d[\xi]}}{\mathcal{G}_d(x)[\xi]} \cdot M_d(x)[\xi]$.
\item\label{l:P_Omega_d ii} If $\mathcal{G}_d(x)[\xi] = 0$, then $y[\xi]$ varies on the set $S_d[\xi]$ defined in (\ref{S_d_xi}).
\end{enumerate}
\end{lemma}
\textbf{Proof.}
By definitions (\ref{Omega_d}), (\ref{M_d}) and (\ref{S_d}) of $\Omega_d$, $M_d$ and $S_d$ respectively, it holds that
\[
    S_d = M_d(\Omega_d).
\]
Then by the unitarity property of $M_d$, we have that
\[
	P_{\Omega_d}(x) = M_d^{-1}\paren{P_{M_d\paren{\Omega_d}}(M_d(x))} = M_d^{-1}\paren{P_{S_d}(M_d(x))}\quad \paren{\forall x\in \Hbb}.
\]
Plugging the formulas of $M_d$, $M_d^{-1}$ and $P_{S_d}$ respectively given by (\ref{M_d}), (\ref{M_d^-1}) and Fact \ref{f:P_S_d} into the above identity, we obtain the characterization of $P_{\Omega_d}$ as claimed.
\hfill$\square$

\begin{lemma}[projection on $\chi$]\label{l:P_chi}
For any $x = \paren{x_{c}}_{c\,\in \mathcal{I}}\in \Hbb$, it holds that
\begin{equation}\label{P_chi}
	P_{\chi}(x) = \left\{\paren{E_{c}\cdot \mathcal{A} \cdot \ee^{{\rm j}\Psi}}_{c\,\in \mathcal{I}} \mid \Psi \in \arg\paren{\sum_{c\,\in \mathcal{I}} \paren{E_{c}\cdot \mathcal{A} \cdot x_{c}}}\right\}.
\end{equation}
\end{lemma}
\textbf{Proof.}
By definition (\ref{chi}) of $\chi$ and the definition of projector in (\ref{def:projector}), $\paren{E_{c}\cdot \mathcal{A} \cdot \ee^{\jj\Psi}}_{c\,\in \mathcal{I}}$ is a projection of $x$ on $\chi$ if and only if $\Psi$ is a solution to the following minimization problem:
\begin{equation}\label{min problem}
	\min_{\Psi \in \mathbb{R}^{n\times n}}\quad  \norm{x - \paren{E_{c}\cdot \mathcal{A} \cdot \ee^{\jj\Psi}}_{c\,\in \mathcal{I}}}^2.
\end{equation}
The objective function of (\ref{min problem}) can be rewritten as
\begin{equation*}
\norm{x - \paren{E_{c}\cdot \mathcal{A} \cdot \ee^{\jj\Psi}}_{c\,\in \mathcal{I}}}^2 = - 2\Re\paren{\ee^{-\jj\Psi} \cdot \sum_{c\,\in \mathcal{I}}\paren{E_{c}\cdot \mathcal{A} \cdot x_{c}}} + C,
\end{equation*}
where $C \equiv \|x\|^2 + \|\paren{E_{c}\cdot \mathcal{A}}_{c\,\in \mathcal{I}}\|^2$ is independent of $\Psi$.
The problem (\ref{min problem}) is hence equivalent to the following one:
\begin{equation}\label{max problem}
	\max_{\Psi \in \mathbb{R}^{n\times n}}\quad  \Re\paren{\ee^{-\jj\Psi} \cdot \sum_{c\,\in \mathcal{I}}\paren{E_{c}\cdot \mathcal{A} \cdot x_{c}}}.
\end{equation}
It is clear that the solution set of the problem (\ref{max problem}) is given by
\[
\left\{\Psi \in \mathbb{R}^{n\times n} \mid \Psi \in \arg\paren{\sum_{c\,\in \mathcal{I}} \paren{E_{c}\cdot \mathcal{A} \cdot x_{c}}}\right\}.
\]
The proof is complete.
\hfill$\square$

\begin{lemma}[projection on $A$]\label{l:P_A}
For any $w = \paren{w_1,w_2,\ldots,w_m}\in \Hbb^m$, it holds that
\begin{equation*}
    P_A(w) = \left[\paren{E_{c}\cdot \overline{z}}_{c \,\in \mathcal{I}}\right]_m\;\mbox{ with }\; \overline{z} \equiv \frac{1}{2}\sum_{c\,\in \mathcal{I}}\paren{E_{c}\cdot \overline{w}_{c}},
\end{equation*}
where $\overline{w} = \paren{\overline{w}_{c}}_{c \,\in \mathcal{I}} \equiv  \paren{1/m}\sum_{d=1}^m w_d$.
\end{lemma}
\textbf{Proof.}
We first note that the set
\begin{equation}\label{L}
L \equiv \left\{[x]_m \in \Hbb^m \mid x \in \Hbb\right\}	
\end{equation}
is a linear subspace of $\Hbb^m$ and contains the set $A$.
By the basic properties of the projection, it holds that
\[
	P_A(w) = P_A\paren{P_L(w)}\quad \forall w = \paren{w_1,w_2,\ldots,w_m} \in\Hbb^m.
\]
By Fact \ref{f:P_D}, the projector $P_L$ admits the following form:
\begin{equation}\label{P_L}
	P_L(w) = [\overline w]_m\; \mbox{ with } \overline{w} \equiv \frac{1}{m}\sum_{d=1}^m w_d.
\end{equation}
This together with the definition of $A$ in (\ref{A&B}) yields that
\begin{equation}\label{key_A}
	P_A(w) = P_A\paren{[\overline w]_m} = \left[P_{\Omega_0}\paren{\overline w}\right]_m.
\end{equation}
The claimed characterization of $P_A$ then follows from (\ref{key_A}) and Lemma \ref{l:P_Omega_0}.
\hfill$\square$

\begin{lemma}[projection on $A_{\chi}$]\label{l:P_A_chi}
For any $w = \paren{w_1,w_2,\ldots,w_m}\in \Hbb^m$, it holds that
\begin{equation*}
	P_{A_\chi}(w) = \left\{\left[\paren{E_{c}\cdot \mathcal{A} \cdot \ee^{{\rm j}\Psi}}_{c \,\in \mathcal{I}}\right]_m \mid \Psi \in \arg\paren{\sum_{c \,\in \mathcal{I}} \paren{E_{c}\cdot \mathcal{A} \cdot \overline{w}_{c}}}\right\}
\end{equation*}
where $\overline{w} = \paren{\overline{w}_{c}}_{c \,\in \mathcal{I}} \equiv  \paren{1/m}\sum_{d=1}^m w_d$.
\end{lemma}
\textbf{Proof.}
The proof is similar to that of Lemma \ref{l:P_A}.
We first observe that the linear subspace $L$ defined in (\ref{L}) contains the set $A_\chi$.
As a consequence, it holds that
\[
	P_{A_\chi}(w) = P_{A_\chi}\paren{P_L(w)} \quad \forall w\in\Hbb^m.
\]
In view of Fact \ref{f:P_D}, the projector $P_L$ admits the explicit form (\ref{P_L}).
This together with the definition of $A_\chi$ in (\ref{A_chi & B_chi}) yields that
\begin{equation}\label{key_A_chi}
	P_{A_\chi}(w) = P_{A_\chi}\paren{[\overline w]_m} = \left[P_\chi\paren{\overline w}\right]_m.
\end{equation}
The claimed characterization of $P_{A_\chi}$ then follows from (\ref{key_A_chi}) and Lemma \ref{l:P_chi}.
\hfill$\square$

\begin{remark}[projections on $B^+$ and $B_{\chi}$]\label{r:P_B^+&P_B_chi}
The projectors $P_{B^+}$ and $P_{B_\chi}$ are analogous to $P_{B}$ in view of Fact \ref{f:P_B}.
\end{remark}

\begin{remark}[nonconvexity]\label{r:nonconvex}
Lemmas \ref{l:P_Omega_d}, \ref{l:P_chi} and \ref{l:P_A_chi} show that the projectors $P_{\Omega_d}$, $P_\chi$ and $P_{A_\chi}$ are not single-valued in general.
This particularly implies that the feasibility models of high-NA phase retrieval formulated in section \ref{subs:FP} are nonconvex.
\end{remark}

\section{Projection algorithms}\label{s:projection algorithm}

Projection methods for phase retrieval can be viewed as descendants of the famous Gerchberg--Saxton algorithm \cite{GerSax72}. 
Its expansions to become the most widely used class of algorithms has been motivated by the rapidly widening scope of phase retrieval applications.
Having calculated the relevant projectors in section \ref{s:prox-operator}, we can implement every projection algorithm for solving the feasibility models formulated in section \ref{subs:FP}.
This section will briefly recall widely known projection methods for solving both two and more-set feasibility problems, typical examples of which are (\ref{PR3}) and (\ref{PR1}), respectively.

Widely known projection methods for solving two-set feasibility are recalled next.

\begin{enumerate}[(i)]
\item The alternating projection (AP) algorithm
\[
T_{\mathtt{AP}}[A,B] \equiv  P_AP_B.
\]
\item The Douglas-Rachford (DR) algorithm
\[
T_{\mathtt{DR}}[A,B] \equiv \frac{1}{2}\paren{R_AR_B+\Id} = P_AR_B - P_B + \Id,
\]
and its Krasnoselski-Mann relaxation (KM-DR algorithm)
\[
T_{\mathtt{KM-DR}} \equiv \beta T_{\mathtt{DR}}[A,B] + (1-\beta)\Id,
\]
where $\beta\in (0,1]$ is the tuning parameter.
\item The \emph{Hybrid Projection-Reflection} (HPR) algorithm \cite[Eq. (19)]{BauComLuk03}:
\[
T_{\mathtt{HPR}} \equiv P_A\left((1+\beta)P_B-\Id\right) - \beta P_B+\Id,
\]
where $\beta\in (0,1]$ is the tuning parameter.
As shown in \cite[Proposition 1]{BauComLuk03}, the HPR algorithm is equivalent to the Fienup's \emph{hybrid input-output} method \cite{Fie82} when $A$ is a linear subspace.
\item The \emph{Relaxed-Averaged-Alternating-Reflections} (RAAR) algorithm \cite{Luk05}:
\begin{eqnarray*}
	T_{\mathtt{RAAR}}[A,B] &\equiv  \frac{\beta}{2}\paren{R_A R_B+\Id} + (1-\beta)
	\\
	&= \beta T_{\mathtt{DR}}[A,B] + (1-\beta)P_B,
\end{eqnarray*}
where $\beta\in (0,1]$ is the tuning parameter.
\item The \emph{Relaxed-Reflect-Reflect} (RRR) algorithm \cite[Algorithm 1]{ElsLanBen18}:
\begin{equation*}
	T_{\mathtt{RRR}}[A,B] \equiv  \beta P_A\left(2P_B-\Id\right) - \beta P_B + \Id,
\end{equation*}
where $\beta\in (0,1]$ is the tuning parameter.
\item The DRAP algorithm \cite{Tha18}:
\[
T_{\mathtt{DRAP}} \equiv P_A\left((1+\beta)P_B-\beta\Id\right)-\beta\left(P_B-\Id\right)
\]
where $\beta\in [0,1]$ is the tuning parameter. This algorithm covers both $T_\mathtt{DR}$ (by setting $\beta = 1$) and $T_\mathtt{AP}$ (by setting $\beta =0$).
When $A$ is affine, $T_{\mathtt{DRAP}}$ is convex combination of these two operators \cite{NguSolVer21}.
The latter also explains its name DRAP which stands for Doughlas-Rachford and Alternating Projection.
\end{enumerate}
\bigskip

Solution algorithms for solving the $(m+1)$-set feasibility are the cyclic projection and the cyclic versions of two-set feasibility based algorithms.

\begin{enumerate}
\item[(vii)] The cyclic projection algorithm
\[
    T_{\mathtt{CP}}[\Omega_0,\Omega_1,\cdots,\Omega_m] \equiv  P_{\Omega_0}P_{\Omega_1}\cdots P_{\Omega_m}.
\]
\item[(viii)] The \emph{cyclic Douglas-Rachford algorithm} proposed and analyzed in the context of convex feasibility \cite{BorTam14}:
\[
	T_{\mathtt{CDR}}[\Omega_0,\Omega_1,\cdots,\Omega_m] \equiv   T_{\mathtt{DR}}[\Omega_0,\Omega_1]
	T_{\mathtt{DR}}[\Omega_1,\Omega_2] \cdots T_{\mathtt{DR}}[\Omega_m,\Omega_0].
\]
\item[(ix)] The \emph{cyclic RAAR algorithm} proposed in the context of low-NA phase retrieval \cite{LukSabTeb19}:
\[
	T_{\mathtt{CRAAR}}[\Omega_0,\Omega_1,\cdots,\Omega_m] \equiv  T_{\mathtt{RAAR}}[\Omega_0,\Omega_1] T_{\mathtt{RAAR}}[\Omega_1,\Omega_2] \cdots T_{\mathtt{RAAR}}[\Omega_m,\Omega_0].
\]
\end{enumerate}

The cyclic projections of $T_{\mathtt{KM-DR}}$, $T_{\mathtt{RRR}}$ and $T_{\mathtt{DRAP}}$ can also be designed similarly.

\begin{remark}[multi-valuedness of projection algorithms]\label{r:not single-valued}
Since the projectors presented in section \ref{s:prox-operator} are potentially multi-valued, the above algorithms built on them are in general not single-valued.
\end{remark}

\begin{remark}[choice of algorithms]\label{r:choose PA} Depending on specific setting of phase retrieval, one algorithm can result in better approximate solutions than another, see also Remark \ref{r:choose FP}.
It is worth mentioning that alternating projection is eventually needed for suppressing noise and model deviations regardless of the chosen algorithm.
\end{remark}

\section{Geometry of high-NA phase retrieval}\label{s:prox-regular}

In this section we analyze the geometry of the high-NA phase retrieval problem.
The sets constituting the feasibility models in section \ref{s:problem formulation} will be shown to be \emph{prox-regular} at the points relevant to our subsequent convergence analysis in section \ref{s:convergence analysis}.
We mention that the prox-regularity property in the context of phase retrieval was first analyzed by Luke \cite[section 3.1]{Luk08}.

\begin{definition}[prox-regularity]\label{d:prox-regular}\cite{PolRocThi00}
A set $\Omega$ is \emph{prox-regular at a point} $\hat{x}\in \Omega$ if the associated projector $P_{\Omega}$ is single-valued around $\hat{x}$.
$\Omega$ is \emph{prox-regular} if it is prox-regular at every of its points.
\end{definition}

\begin{example}[prox-regularity of $\Omega_0$, $A$ and $D$]\label{e:prox-regular}
Any closed and convex set is prox-regular \cite{VA}. 
In particular, the linear subspaces $\Omega_0$, $A$ and $D$ defined in (\ref{Omega_0}) and (\ref{A&B}) are prox-regular.
\end{example}

The next two assertions follow from the definition of prox-regularity. Recall the notation $[\,\cdot\,]_p$ in (\ref{duplication}).

\begin{fact}\label{f:prox-regular of duplication} Let $\Omega \subset \Hcal$ be prox-regular at a point $\hat{x}\in \Omega$ and $p\ge 2$ be an integer.
Then the set $\mathbf{\Omega}\equiv  \parennn{[x]_p \in \Hcal^p \mid x \in \Omega}$ is prox-regular at $[\hat{x}]_p$.
\end{fact}
\textbf{Proof.}
By Definition \ref{d:prox-regular}, there is a neighborhood $U$ of $\hat x$ on which $P_\Omega$ is single-valued.
Let us define the set $\mathcal{U}\subset \Hcal^p$ by
\begin{equation}\label{U for duplication}
	\mathcal{U} \equiv  \parennn{[\hat x]_p + (r_1,r_2,\ldots,r_p) \in \Hcal^p\; \mid\; \hat x + \frac{1}{p}\sum_{k=1}^pr_k \in U}.
\end{equation}
Note that $\mathcal{U}$ is a neighborhood of $[\hat{x}]_p$ since $U$ is a neighborhood of $\hat x$.
It suffices to check that $P_{\mathbf{\Omega}}$ is single-valued on $\mathcal{U}$.
Indeed, take an arbitrary point 
\[
(x_1,x_2,\ldots,x_p) = [\hat x]_p + (r_1,r_2,\ldots,r_p) \in \mathcal{U}.
\]
Then in view of (\ref{U for duplication}), it holds that
\[
\bar{x} \equiv \frac{1}{p}\sum_{k=1}^p x_k = \hat x + \frac{1}{p}\sum_{k=1}^pr_k \in U,
\]
and hence $P_\Omega(\bar{x})$ is singleton since $P_\Omega$ is single-valued on $U$.
Using the reasoning in the proof of Lemma \ref{l:P_A}, we have $P_{\mathbf{\Omega}}(x_1,x_2,\ldots,x_p) = \parennn{\parenn{P_\Omega(\bar{x})}_p}$ which is singleton.
Hence $P_{\mathbf{\Omega}}$ is singled-valued on $\mathcal{U}$ and the proof is complete.
\hfill$\square$

\begin{fact}[prox-regularity of products]\label{f:prox-regular of product} For each $k=1,2,\ldots,p$ let $\Omega_k$ be prox-regular at $\hat{x}_k$.
Then the product set $\prod_{k=1}^p \Omega_k$ is prox-regular at $(\hat{x}_1,\hat{x}_2,\ldots,\hat{x}_p)$. 
\end{fact}
\textbf{Proof.}
The proof follows from the definition of prox-regularity and the separation property of projection on product sets.
\hfill$\square$

We can now analyze the prox-regularity of the other sets defined in section \ref{subs:FP}.

\begin{proposition}[prox-regularity of $\Omega_d$]\label{p:prox-regular of Omega_d}
For each $d=1,2,\ldots,m$ the set $\Omega_d$ defined in (\ref{Omega_d}) is prox-regular at every point $\hat x \in \Omega_d$ with $\mathcal{G}_d(\hat x)$ nonzero everywhere.
\end{proposition}
\textbf{Proof.}
Consider a point $\hat x \in \Omega_d$ with $\mathcal{G}_d(\hat x)$ nonzero everywhere.
By Definition \ref{d:prox-regular}, it suffices to find a neighborhood of $\hat x$ on which $P_{\Omega_d}$ is single-valued.
Let us define the set $U_d$ by
\begin{equation}\label{U for Omega_d}
	U_d \equiv  \parennn{\hat x + r \mid r \in \mathcal{H},\, \mathcal{G}_d(r) < \mathcal{G}_d(\hat{x})}.
\end{equation}
Since $\mathcal{G}_d$ is continuous by Fact \ref{f:continuous of G_d}, it holds that
\begin{equation*}
	\mathcal{G}_d(r) \to 0 \mbox{ in } \mathbb{R}_+^{n\times n} \quad \mbox{as}\quad r \to 0 \mbox{ in } \mathcal{H}.
\end{equation*}
This together with $\mathcal{G}_d(\hat{x})$ being nonzero everywhere implies that $U_d$ is a neighborhood of $\hat x$.
We will show that $P_{\Omega_d}$ is single-valued on $U_d$.
Indeed, let us take an arbitrary point $x = \hat x + r \in U_d$ and first check that $\mathcal{G}_d(x) \neq 0$ for all entries.
Using (\ref{G_d}), the linearity of $M_d$ and the triangle inequality successively, we have that
\begin{eqnarray}\nonumber	
	\mathcal{G}_d(x) = \mathcal{G}_d(\hat x+r) 
	&= \sqrt{\sum_{c\,\in \mathcal{I}}\left|M_d\paren{\hat{x}+r}_{c}\right|^2}
	\\\label{>0 for Omega_d}
	&= \sqrt{\sum_{c\,\in \mathcal{I}}\left|M_d\paren{\hat{x}}_{c} + M_d\paren{r}_{c}\right|^2}
	\\\nonumber
	&\ge \sqrt{\sum_{c\,\in \mathcal{I}}\paren{ \left|M_d\paren{\hat{x}}_{c}\right|
	- \left|M_d\paren{r}_{c}\right|}^2}.
\end{eqnarray}
Suppose on the contrary that $\mathcal{G}_d(x)[\xi] = 0$ for some index $\xi\in \mathcal{J}$.
Then (\ref{>0 for Omega_d}) implies that
\[
\left|M_d\paren{\hat{x}}_{c}\right|[\xi] = \left|M_d\paren{r}_{c}\right|[\xi]\quad \forall c \in \mathcal{I}.
\]
This in particular yields $\mathcal{G}_d(\hat x)[\xi] = \mathcal{G}_d(r)[\xi]$ which is a contradiction to (\ref{U for Omega_d}) as $\hat x+r\in U_d$.
Hence we have $\mathcal{G}_d(x) \neq 0$ for all entries as claimed.
Now by Lemma \ref{l:P_Omega_d}, $P_{\Omega_d}(x)$ is the singleton
$\parennn{M_d^{-1}\paren{y}}$, where $y$ is uniquely determined.
The proof is complete.
\hfill$\square$

\begin{proposition}[prox-regularity of $\chi$]\label{p:prox-regular of chi}
Suppose that the amplitude $\mathcal{A}$ is nonzero everywhere.
Then the set $\chi$ defined in (\ref{chi}) is prox-regular.
\end{proposition}
\textbf{Proof.}
Let us consider an arbitrary point $\hat x = \paren{\hat{x}_{c}}_{c\,\in \mathcal{I}} \in \chi$.
By Definition \ref{d:prox-regular}, it suffices to find a neighborhood of $\hat x$ on which $P_\chi$ is single-valued.
Let us define the set $U_\chi$ by
\begin{equation}\label{U for chi}
U_\chi \equiv  \parennn{\hat x + \paren{r_{c}}_{c\,\in \mathcal{I}} \in \mathcal{H}\; \mid\; \sqrt{\sum_{c\,\in \mathcal{I}} \parennnn{r_{c}}^2} < \sqrt{2}\, \mathcal{A}}.
\end{equation}
Since $\mathcal{A}$ is nonzero everywhere, the set $U_\chi$ defined in (\ref{U for chi}) is a neighborhood of $\hat x$.
We will show that $P_\chi$ is single-valued on $U_\chi$.
Take an arbitrary point $x = \paren{x_{c}}_{c\,\in \mathcal{I}} = \hat x + \paren{r_{c}}_{c\,\in \mathcal{I}} \in U_\chi$.
Then using the triangle inequality, the Cauchy-Schwartz inequality, (\ref{chi}), (\ref{sum=2}) and (\ref{U for chi}) successively, we get that
\begin{eqnarray*}	
&\parennnn{\sum_{c\,\in \mathcal{I}}\paren{E_{c}\cdot \mathcal{A} \cdot x_{c}}} 
= \parennnn{\sum_{c\,\in \mathcal{I}}\paren{E_{c}\cdot \mathcal{A} \cdot (\hat{x}_{c} + r_{c})}}
\\
&\ge \parennnn{\sum_{c\,\in \mathcal{I}} \paren{E_{c} \cdot \mathcal{A} \cdot \hat{x}_{c}}} - 
\parennnn{\sum_{c\,\in \mathcal{I}}\paren{E_{c} \cdot \mathcal{A}\cdot r_{c}}}
\\
&\ge \parennnn{\sum_{c\,\in \mathcal{I}} \paren{E_{c} \cdot \mathcal{A} \cdot \hat{x}_{c}}} - \sqrt{\sum_{c\,\in \mathcal{I}} \paren{E_{c} \cdot \mathcal{A}}^2\,}\, \cdot \sqrt{\sum_{c\,\in \mathcal{I}} |r_{c}|^2}\label{>0 for chi}
\\
&= \sum_{c\,\in \mathcal{I}} \paren{E_{c} \cdot \mathcal{A}}^2 
- \sqrt{\sum_{c\,\in \mathcal{I}} \paren{E_{c} \cdot \mathcal{A}}^2\,}\, \cdot \sqrt{\sum_{c\,\in \mathcal{I}} |r_{c}|^2}
\\
&= 2 \mathcal{A}^2 
- \sqrt{2}\; \mathcal{A} \cdot \sqrt{\sum_{c\,\in \mathcal{I}} |r_{c}|^2}\; > 0.
\end{eqnarray*}
This implies that $\sum_{c\,\in \mathcal{I}} \paren{E_{c}\cdot \mathcal{A} \cdot x_{c}}$ is nonzero everywhere.
Hence, by Lemma \ref{l:P_chi}, $P_{\chi}(x)$ is the singleton $\parennn{\paren{E_{c}\cdot \mathcal{A}\cdot \ee^{\jj\Psi}}_{c\,\in\mathcal{I}}}$, where $\Psi\in \mathbb{R}^{n\times n}$ is uniquely given by (\ref{P_chi}).
The proof is complete.
\hfill$\square$

\begin{proposition}[prox-regularity of $A_\chi$]\label{p:prox-regular of A_chi}
Suppose that the amplitude $\mathcal{A}$ is nonzero everywhere.
Then the set $A_\chi$ defined in (\ref{A_chi & B_chi}) is prox-regular at every point $\parenn{\hat{x}}_m$ with $\hat{x} \in \chi$.
\end{proposition}
\textbf{Proof.}
The proof follows from Proposition \ref{p:prox-regular of chi} and Fact \ref{f:prox-regular of duplication}.
\hfill$\square$

\begin{proposition}[prox-regularity of $B$, $B^+$ and $B_\chi$]\label{p:prox-regular of B, B^+, B_chi}
The following statements hold true.
\begin{enumerate}[(i)]
\item\label{prox-regular of B} The set $B$ defined in (\ref{A&B}) is prox-regular at every point $\paren{\hat{x}_1,\hat{x}_2,\ldots,\hat{x}_m} \in B$ with $\mathcal{G}_d(\hat{x}_d)$ nonzero everywhere $(\forall d=1,2,\ldots,m)$.
\item\label{prox-regular of B^+} The set $B^+$ defined in (\ref{A&B}) is prox-regular at every point $\paren{\hat{x},\hat{x}_1,\hat{x}_2,\ldots,\hat{x}_m} \in B^+$ with $\mathcal{G}_d(\hat{x}_d)$ nonzero everywhere $(\forall d=1,2,\ldots,m)$.
\item\label{prox-regular of B_chi} Suppose that the amplitude $\mathcal{A}$ is nonzero everywhere.
Then the set $B_\chi$ defined in (\ref{A_chi & B_chi}) is prox-regular at every point $\paren{\hat{x},\hat{x}_1,\hat{x}_2,\ldots,\hat{x}_m} \in B_\chi$ with $\mathcal{G}_d(\hat{x}_d)$ nonzero everywhere $(\forall d=1,2,\ldots,m)$.
\end{enumerate}
\end{proposition}
\textbf{Proof.}
\ref{prox-regular of B} By Proposition \ref{p:prox-regular of Omega_d}, for each $d=1,2,\ldots,m$ there exists a neighborhood $U_d$ of $\hat{x}_d$ on which $P_{\Omega_d}$ is single-valued.
This combined with Fact \ref{f:P_B} yields that $P_B= \prod_{d=1}^m P_{\Omega_d}$ is single-valued in the neighborhood $\prod_{d=1}^m U_d$ of $\paren{\hat{x}_1,\hat{x}_2,\ldots,\hat{x}_m}$.
This yields the prox-regularity of $B$ at this point as claimed.
\ref{prox-regular of B^+} This part is also encompassed by part \ref{prox-regular of B} since $\Omega_0$ is prox-regular in view of Example \ref{e:prox-regular}.
\ref{prox-regular of B_chi} Thanks to Proposition \ref{p:prox-regular of chi}, there exists a neighborhood $U_\chi$ of $\hat{x}$ on which $P_\chi$ is single-valued.
By Proposition \ref{p:prox-regular of Omega_d}, for each $d=1,2,\ldots,m$ there exists a neighborhood $U_d$ of $\hat{x}_d$ on which $P_{\Omega_d}$ is single-valued.
We thus have in view of Remark \ref{r:P_B^+&P_B_chi} that $P_{B_\chi}= P_\chi \times \prod_{d=1}^m P_{\Omega_d}$ is single-valued on the neighborhood $U_\chi \times \prod_{d=1}^m U_d$ of $\paren{\hat{x}, \hat{x}_1, \hat{x}_2, \ldots, \hat{x}_m}$.
This yields the prox-regularity of $B_\chi$ at this point as claimed.
The proof is complete.
\hfill$\square$

\begin{remark}\label{r:in aperture for Omega_0}
The condition that $\mathcal{A}$ are nonzero everywhere imposed in Propositions \ref{p:prox-regular of chi}, \ref{p:prox-regular of A_chi} and \ref{p:prox-regular of B, B^+, B_chi}\ref{prox-regular of B_chi} physically means that the entire aperture of the imaging system is illuminated.
\end{remark}

\section{Convergence analysis}\label{s:convergence analysis}

The feasibility models of high-NA phase retrieval formulated in section \ref{subs:FP} are nonconvex (Remark \ref{r:nonconvex}) and hence the projection algorithms are not Fej\'er monotone, indeed not even single-valued (Remark \ref{r:not single-valued}).
As a result, tools in convex analysis and monotone operator theory (for example, \cite{BauBor96,BauCom17}) are not applicable to the problem under study.
In this paper, we follow the analysis scheme of \cite{LukNguTam18} according to which convergence of iterative sequences generated by a fixed point operator $T$ is guaranteed by the \emph{pointwise almost averagedness} of $T$ and the \emph{metric subregularity} of the mapping $\Id-T$ on the relevant regions.
The contribution of this section concerns the first condition of convergence.
The almost averagedness property of projection algorithms will be derived from the geometry of the high-NA phase retrieval problem analysed in section \ref{s:prox-regular}.

Although being derived from the general scheme of \cite{LukNguTam18}, convergence analysis is different for each projection method, depending on its fixed point set and its complexity, especially for solving nonconvex and inconsistent feasibility problems.
In this section, we analyze the alternating projection algorithm for high-NA phase retrieval in the presence of noise.
We consider the two-set feasibility model (\ref{PR3}) in the \emph{inconsistent} setting, \textit{i.e.}, the sets do not intersect.
It is worth emphasizing that the class of projection methods for high-NA phase retrieval is first considered in this paper, and thus the obtained results are new from the application point of view, even in the consistent case.

\subsection{Pointwise almost averagedness}\label{subs:a.a.}

The following property is taken from Definition 2.2 and Proposition 2.1 of \cite{LukNguTam18}.

\begin{definition}[pointwise almost averaged mappings]\label{d:p.a.a.}
A fixed point mapping $T:\Hcal \rightrightarrows \Hcal$ is \emph{pointwise almost averaged} at a point $y$ on a set $\Omega \subset \Hcal$ with violation $\varepsilon\ge 0$ and averaging constant $\alpha \in (0,1)$ if for all $y^+\in T(y)$, $z\in \Omega$ and $z^+\in T(z)$, it holds that
\begin{equation*}
	\norm{z^+- y^+}^2 \le \paren{1+\varepsilon}\norm{z-y}^2 - \frac{1-\alpha}{\alpha}\norm{(z^+-z)-(y^+-y)}^2.
\end{equation*}
\end{definition}

When the violation $\varepsilon=0$, the quantifiers `almost' and `violation' in Definition \ref{d:p.a.a.} are dropped and the property goes back to the conventional averagedness property, see, for example, \cite{BauCom17}.
When the property holds for every point $y\in \Omega$ with the same violation and averaging constant, the quantifiers `pointwise' and `at a point' in Definition \ref{d:p.a.a.} are dropped.
The property is well defined for any averaging constant $\alpha>0$, not necessarily limited to $\alpha\in (0,1)$ though the latter is often of the main interest.

\begin{example}[projection on convex sets]\label{e:project on convex set}
The projectors associated with closed and convex sets are globally averaged with averaging constant $\alpha=1/2$ (\textit{i.e.}, firmly nonexpansive), see, for example, \cite[Theorem 2.2.21]{Ceg12}.
\end{example}

The following statement is a consequence of widely known results concerning projections on nonconvex sets, see, for example, \cite[Theorem 2.14]{HesLuk13}.

\begin{proposition}[projection on prox-regular sets]\label{p:project on prox-regular set}
Let $\Omega$ be closed and prox-regular at $\hat x \in \Omega$.
Then given an arbitrarily small number $\varepsilon>0$, there exists a neighborhood of $\hat{x}$ (depending on $\varepsilon$) on which $P_\Omega$ is almost averaged with violation $\varepsilon$ and averaging constant $\alpha=1/2$.
\end{proposition}

The next property of pointwise almost averaged mappings is needed \cite[Proposition 2.4(ii)]{LukNguTam18}.
The version specialized to the problem (\ref{PR3}) is presented here for brevity.

\begin{proposition}[pointwise almost averagedness of composite mappings] \label{p:composition of a.a.}
Let\\
$\mmap{T_k}{\Hcal}{\Hcal}$ for $k=1,2$ be pointwise almost averaged on $U_k$ at all $y_k\in S_k$ with violation $\varepsilon_k\ge 0$ and averaging constant $\alpha_k\in (0,1)$.
If $T_2\paren{U_2} \subseteq U_1$ and $T_2\paren{S_2} \subseteq S_1$, then the composite mapping $T\equiv T_1\circ T_{2}$ is pointwise almost averaged on $U_2$ at all $y\in S_2$ with violation $\varepsilon$ and averaging constant $\alpha$ given by
\begin{equation*}
	\varepsilon = \varepsilon_1 + \varepsilon_2 + \varepsilon_1\varepsilon_2;\quad 
	\alpha = \frac{2\max\{\alpha_1,\alpha_2\}}{1 + \max\{\alpha_1,\alpha_2\}}.
\end{equation*}
\end{proposition}

The next result links the prox-regularity of the sets in (\ref{PR3}) with the almost averagedness of the alternating projection operator.

\begin{proposition}[almost averagedness of $T_{\mathtt{AP}}$]\label{p:a.a. of AP} Let $\widehat{\mathbf{b}} \equiv \paren{\hat{x}_1,\hat{x}_2,\ldots,\hat{x}_m} \in B$, where $\hat{x}_d \in \Omega_d$ with $\mathcal{G}_d(\hat{x}_d)$ nonzero everywhere $(\forall d=1,2,\ldots,m)$.
Then given any number $\varepsilon>0$, there is a neighborhood of $\widehat{\mathbf{b}}$, denoted by $U_{\varepsilon}(\widehat{\mathbf{b}})$, on which the alternating projection operator $T_{\mathtt{AP}}\equiv P_AP_B$ associated with (\ref{PR3}) is almost averaged with violation $\varepsilon$ and averaging constant $\alpha=2/3$.
\end{proposition}
\textbf{Proof.}
By Proposition \ref{p:prox-regular of Omega_d}, the sets $\Omega_d$ are prox-regular at $\hat{x}_d$ as $\mathcal{G}_d(\hat{x}_d)$ are nonzero everywhere ($\forall d=1,2,\ldots,m$).
Thanks to Fact \ref{f:prox-regular of product}, the set $B$ is prox-regular at $\widehat{\mathbf{b}}$.
Then by Proposition \ref{p:project on prox-regular set}, there exists a neighborhood $U_{\varepsilon}(\widehat{\mathbf{b}})$ of $\widehat{\mathbf{b}}$ on which the projector $P_B$ is almost averaged with violation $\varepsilon$ and averaging constant $1/2$.
On the other hand, since $A$ is convex in view of Example \ref{e:prox-regular}, the projector $P_A$ is globally averaged with averaging constant $1/2$ (\textit{i.e.}, firmly nonexpansive) in view of Example \ref{e:project on convex set}.
Thus by Proposition \ref{p:composition of a.a.}, the composite mapping $T_{\mathtt{AP}}\equiv P_AP_B$ is almost averaged on $U_{\varepsilon}(\widehat{\mathbf{b}})$ with violation $\varepsilon$ and averaging constant $\alpha=2/3$ as claimed.
\hfill$\square$

\subsection{Convergence statements}\label{subs:convergence results}

The goal of this section is to combine the results of section \ref{subs:a.a.} with the analysis scheme of \cite[section 2.2]{LukNguTam18} to obtain convergence criteria for the alternating projection algorithm for solving (\ref{PR3}) in the inconsistent setting.
The following notion of metric subregularity is a cornerstone of variational analysis and optimization theory with many important applications, such as in establishing calculus rules for subdifferentials and coderivatives \cite{Iof17,Mor06,VA} and in analyzing stability and convergence of numerical algorithms, see, for example, \cite{DonRoc14,KlaKum02}.

\begin{definition}[metric subregularity on a set]
A set-valued mapping $\Theta:\mathcal{H}\rightrightarrows\mathcal{H}$ is \emph{metrically subregular} on $U \subset \mathcal{H}$ for $\hat{y} \in \mathcal{H}$ relative to $\Lambda\subset \mathcal{H}$ with \emph{modulus} $\kappa>0$ if
\[
\kappa\, \dist(x,\Theta^{-1}(\hat{y})\cap \Lambda) \le \dist(\hat{y},\Theta(x)), \quad \forall x\in U\cap \Lambda.
\]
When $U$ is some neighborhood of a point $\hat{x}\in \Theta^{-1}(\hat{y})$, the property is called \emph{metric subregularity of $\Theta$ at $\hat{x}$ for $\hat{y}$ relative to $\Lambda$}.
\end{definition}

The next lemma is a specification of \cite[Corollary 2.3]{LukNguTam18} to our target application.

\begin{lemma}[linear convergence with metric subregularity]\label{l:abstract convergence} Let $\mmap{T}{\Hcal}{\Hcal}$ be a fixed point operator with $\Fix T$ closed, $\Lambda\subset \mathcal{H}$ with $T(\Lambda) \subset \Lambda$, $\hat{x} \in \Lambda \cap \Fix T$ and $U$ a neighborhood of $\hat{x}$ with $T(U)\subset U$.
Suppose that
\begin{enumerate}[(a)]
\item\label{t:Tconv a} $T$ is pointwise almost averaged at $\hat{x}$ on $\Lambda \cap U$ with violation $\varepsilon\ge 0$ and averaging constant $\alpha\in (0,1)$;
\item\label{t:Tconv b} the mapping $\Id - T$ is metrically subregular on $U$ for $0$ relative to $\Lambda$ with modulus $\kappa>\sqrt{\varepsilon\alpha/(1-\alpha)}$.
\end{enumerate}
Then every iterative sequence generated by $T$ with the initial point in $\Lambda \cap U$ converges linearly to a point in $\Fix T$ with rate at most (worst) $c \equiv \sqrt{1+\varepsilon-\kappa^2(1-\alpha)/\alpha}\;<1$.
\end{lemma}

We are now ready to formulate the main convergence results.

\begin{theorem}[linear convergence of $T_{\mathtt{AP}}$ for (\ref{PR3})]\label{t:convergence inconsistent}
Let $\widehat{\mathbf{a}} \in A$ be a fixed point of $T_{\mathtt{AP}} \equiv P_AP_B$ and suppose that $P_B(\widehat{\mathbf{a}}) = \parennn{\widehat{\mathbf{b}} \equiv \paren{\hat{x}_1,\hat{x}_2,\ldots,\hat{x}_m}}$ is singleton with $\mathcal{G}_d(\hat{x}_d)$ nonzero everywhere $(\forall d=1,2,\ldots,m)$.
Given a number $\varepsilon>0$, let $U_{\varepsilon}(\widehat{\mathbf{b}})$ be the neighborhood of $\widehat{\mathbf{b}}$ on which $T_{\mathtt{AP}}$ is almost averaged with violation $\varepsilon$ and averaging constant $\alpha=2/3$ as determined by Proposition \ref{p:a.a. of AP}.
Suppose further that $\widehat{\mathbf{a}} \in U_{\varepsilon}(\widehat{\mathbf{b}})$, $T_{\mathtt{AP}}(A\cap U_{\varepsilon}(\widehat{\mathbf{b}})) \subset U_{\varepsilon}(\widehat{\mathbf{b}})$ and the mapping $\Theta \equiv \Id - T_{\mathtt{AP}}$ is metrically subregular on $U_{\varepsilon}(\widehat{\mathbf{b}})$ for $0$ relative to $A$ with modulus $\kappa>\sqrt{2\varepsilon}$.
Then every iterative sequence generated by $T_{\mathtt{AP}}$ with the initial point in $A\cap U_{\varepsilon}(\widehat{\mathbf{b}})$ converges linearly to a point in $\Fix T_{\mathtt{AP}}$ with rate at most $c \equiv \sqrt{1+\varepsilon-\kappa^2/2\,}\;<1$.
\end{theorem}
\textbf{Proof.}
The assumption $\widehat{\mathbf{a}}\in U_{\varepsilon}(\widehat{\mathbf{b}})$ ensures that $T_{\mathtt{AP}}$ is pointwise almost averaged at $\widehat{\mathbf{a}}$ on $U_{\varepsilon}(\widehat{\mathbf{b}})$ with violation $\varepsilon$ and averaging constant $\alpha=2/3$ in view of Proposition \ref{p:a.a. of AP}.
Hence all the assumptions of Lemma \ref{l:abstract convergence} are satisfied with $\Lambda=A$ and $U = U_{\varepsilon}(\widehat{\mathbf{b}})$ and the convergence statement follows as claimed.
\hfill$\square$

We next explain and remark on the assumptions imposed in Theorem \ref{t:convergence inconsistent}.

\begin{remark}\label{r:technical assumptions}
It is important to keep in mind that $U_{\varepsilon}(\widehat{\mathbf{b}})$ is not limited to some ball centered at $\widehat{\mathbf{b}}$.
It can be an unbounded set, see, for example, the typical intuitive example of phase retrieval in \cite[Figure 3 and Example 3.9(ii)]{LukMar20}.
This in particular makes the assumption $\widehat{\mathbf{a}} \in U_{\varepsilon}(\widehat{\mathbf{b}})$ not restrictive.
A more general notion than prox-regularity called \emph{regularity at a distance} was proposed in \cite{LukMar20} for analyzing the RAAR algorithm for nonconvex and inconsistent feasibility.
However, we are unable to verify that property for the high-NA phase retrieval problem and hence do not apply it to the analysis in this application paper to avoid further unverifiable assumptions.
\end{remark}

\begin{remark} Since the set $B$ in (\ref{PR3}) is compact, every iterative sequence generated by the alternating projection methods has a subsequence converging to a point in $\Fix T_{\mathtt{AP}}$, a local best approximation point to $B$.
Theorem \ref{t:convergence inconsistent} provides sufficient conditions for local linear convergence of the algorithm around a single fixed point. 
Its assumptions can be strengthened for all fixed points of $T_{\mathtt{AP}}$ to yield global convergence of the algorithm, but the quality of the fixed point it converges to and the convergence rate in general depend on where it starts as the problem is nonconvex.
However, such additional assumptions would be unverifiable for the high-NA phase retrieval problem, we chose not to include them in this application paper.
\end{remark}

\begin{remark}[necessity of metric subregularity]
As mentioned early  this section there are two groups of properties often required to prove convergence of nonconvex optimization algorithms.
The geometry of the high-NA phase retrieval problem analyzed in section \ref{s:prox-regular} yields the first one -- \emph{pointwise almost averagedness}.
It has been known that the second one -- \emph{metric subregularity} is difficult to verify, but as been shown in \cite{LukTebTha20} this condition is not only \emph{sufficient} but also {\em necessary} for local linear convergence.
\end{remark}

The mathematical complication of Theorem \ref{t:convergence inconsistent} is mainly due to the \emph{inconsistency} of the problem under study.
In the consistent setting, it reduces to the following much simpler form, where the metric subregularity of $\Id-T_{\mathtt{AP}}$ also reduces to the more intuitive notion called \emph{subtransversality} of the collection of sets $\{A,B\}$ at the intersection point.
For cartoon model of phase retrieval consisting of two (products of) spheres, the subtransversality property is satisfied except when they are tangent.
The proof of the next statement follows from the one of Theorem \ref{t:convergence inconsistent} and is left for brevity.

\begin{corollary}[linear convergence of $T_{\mathtt{AP}}$ for consistent (\ref{PR3})]\label{c:convergence consistent}
Consider the problem (\ref{PR3}) with $A \cap B \neq \emptyset$.
Let $\widehat{\mathbf{a}} \equiv [\hat{x}]_m \in A \cap B$ with $\mathcal{G}_d(\hat{x})$ nonzero everywhere $(\forall d=1,2,\ldots,m)$.
Given a number $\varepsilon>0$, let $\mathbb{B}_{\varepsilon}(\widehat{\mathbf{a}})$ be the ball on which $T_{\mathtt{AP}}$ is almost averaged with violation $\varepsilon$ and averaging constant $\alpha=2/3$ as determined by Proposition \ref{p:a.a. of AP}.
Suppose that the mapping $\Theta \equiv \Id - T_{\mathtt{AP}}$ is metrically subregular at $\widehat{\mathbf{a}}$ for $0$ relative to $A$ with modulus $\kappa>\sqrt{2\varepsilon}$.
Then every iterative sequence generated by $T_{\mathtt{AP}}$ with the initial point in $A\cap \mathbb{B}_{\varepsilon}(\widehat{\mathbf{a}})$ converges linearly to a point in $\Fix T_{\mathtt{AP}}$ with rate at most $c \equiv \sqrt{1+\varepsilon-\kappa^2/2\,}\;<1$.
\end{corollary}

It is worth mentioning that the technical assumptions imposed in Theorem \ref{t:convergence inconsistent} and Corollary \ref{c:convergence consistent} concerning high-NA phase retrieval also remain unverifiable for the low-NA problem.

\section{Numerical simulations}\label{s:numerical result}

The goal of this section is to demonstrate that the new mathematical analysis obtained in this paper enables us to apply the class of projection algorithms for solving the high-NA phase retrieval problem.
In contrast to a vast number of existing solution methods for low-NA phase retrieval, very few algorithms have been proposed for the high-NA case.
The \textit{Vectorial PSF model-based Alternating Minimization} (VAM) algorithm was proposed in \cite{NguSolVer20}.
It outperforms several available high-NA phase retrieval approaches, including the \textit{Scalar PSF model-based Alternating Minimization} (SAM) algorithm of Hanser \textit{et al.} \cite{HanGusAgaSed03} which is limited in accuracy due to model deviations, and the modal-based approach through the use of extended Nijboer–Zernike expansion of Braat \textit{et al.} \cite{BraDirJanHavNes05} which is of high computational complexity and excludes applications with discontinuous phase.
The VAM algorithm is nothing else, but the alternating projection method applied to the feasibility model (\ref{PR3}).
The projectors computed in section \ref{s:prox-operator} enable the implementation of every projection method (not only those mentioned in section \ref{s:projection algorithm}) for solving every corresponding feasibility model formulated in section \ref{s:problem formulation}.
This section aims at demonstrating the improved performance of more delicate projection algorithms over available solution methods for high-NA phase retrieval.
As projection methods have not been applied to high-NA phase retrieval before, their comparison is not a goal of this paper, which instead establishes groundwork enabling the implementation and analysis of this efficient class of solution methods for high-NA phase retrieval.

\begin{table}[b!]
\caption{\label{tbl:parameters}Parameters used in numerical simulations: NA -- numerical aperture, $\lambda$ -- wavelength of illumination light ($\mu$m), $s$ -- pixel size ($\mu$m), $m$ -- number of images, $n \times n$ -- image size (pixels), $w$ -- noise model, and SNR -- signal-to-noise ratio (decibels).}
\vspace{.25cm}
\begin{indented}
\item[]
\begin{tabular}[1\baselineskip]{r|ccccccc}
\textbf{Parameter} & NA & $\lambda$ & $s$ & $m$ & $n \times n$ & $w$ & SNR\cr 
\hline
\textbf{Value} & 0.95 & 0.3 & 0.06 & $7$ & $128 \times 128$ & Gaussian & 30 dB\\
\end{tabular}
\end{indented}
\end{table}

\begin{table}[b!]
\caption{\label{tbl:beta and iter}The number of iterations (the second row) and the parameter $\beta$ (the third row) of the algorithms used in numerical simulations. The averaged RMS errors of phase retrieval over 75 phase realizations are presented in the last row.}
\vspace{.25cm}
\begin{indented}
\item[]
\begin{tabular}[1\baselineskip]{r|ccccccc}
\textbf{Algorithm} & SAM & VAM & DRAP & RAAR & VAM$_+$ & DRAP$_+$ & RAAR$_+$\\
\hline
\textbf{\#Iterations} & 100 & 100 & 30+20 & 30+20 & 100 & 30+20 & 30+20\\
\textbf{Parameter $\beta$} & \  & \  & 0.95 & 0.95 & \  & 0.95 & 0.95\\
\textbf{Error (\%)} & 8.47 & 7.69 & 6.14 & 5.98 & 6.82 & 4.68 & 4.69\\
\end{tabular}
\end{indented}
\end{table}

We consider the practically relevant simulation setting of high-NA phase retrieval as in \cite[section 5]{NguSolVer20} where the vectorial PSF (\ref{vector PSF}) is taken as the forward imaging model for generating the images.
The simulated imaging system has circular aperture with the amplitude $\mathcal{A}$ being the two-dimensional Gaussian distribution truncated at $0.5$ on the boundary.
We do 75 experiments for different phase realizations with values in $[-\pi,\pi]$.
Each data set consists of seven out-of-focus PSF images which are uniformly separated by one depth of focus along the optical axis.
A schematic diagram of this phase retrieval setup can be seen, for example, in \cite[Figure 1]{NguSolVer20}.
The generated PSF images after being normalized to unity energy are corrupted by additive white Gaussian noise with signal-to-noise ratio (SNR) $30$ decibels (dB).
Recall that $\mbox{SNR} = 10\ln\left({P}/{P_0}\right)$, where $P$ and $P_0$ are the powers of the signal and the noise, respectively.
The parameters used in the simulation experiments are summarized in Table \ref{tbl:parameters}.
The quality of phase retrieval is  measured by the relative Root Mean Square (RMS) error ${\|\widehat{\Phi}-\Phi\|}\big{/}{\norm{\Phi}}$, where $\Phi$ and $\widehat{\Phi}$ are the simulation and the retrieved phase aberrations, respectively.
As phase retrieval is ambiguous up to at least a global phase shift (a piston term or the first Zernike mode), the norms of the phases are computed with the piston terms removed.

\begin{figure}[b!]
\centering
\includegraphics[scale=.8]{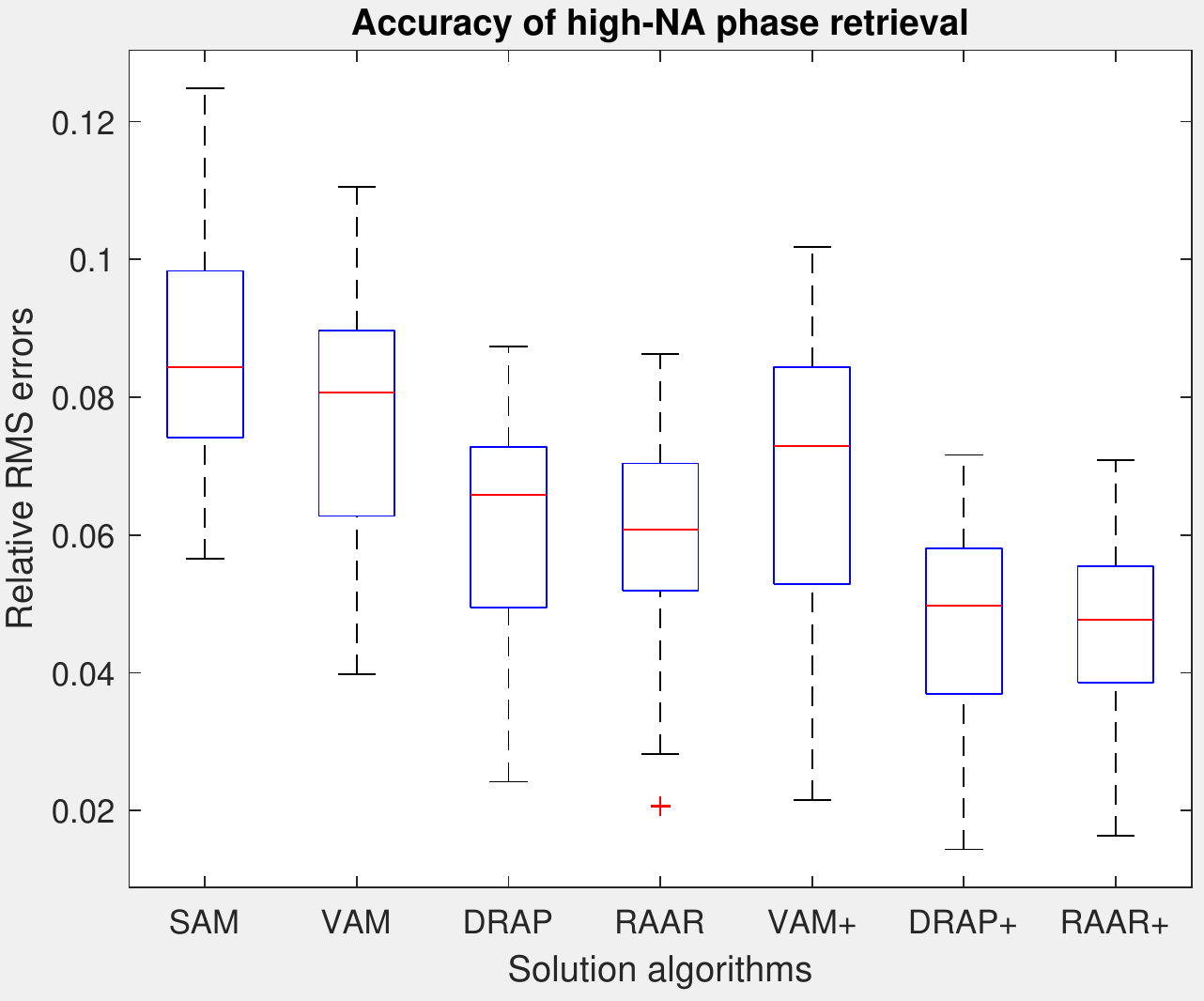}
\caption{The box-plots show the improved performance of the RAAR and DRAP algorithms over available high-NA phase retrieval methods, including SAM \cite{HanGusAgaSed03} and VAM \cite{NguSolVer20}. Each box-plot summarizes the numerical results in relative RMS errors of seventy-five examples with different phase realizations taking values in $[-\pi,\pi]$. The RAAR algorithm yields phase retrieval with the smallest RMS error on average, $5.98$\% compared to $8.47$\% of SAM, $7.69$\% of VAM and $6.14$\% of DRAP. RAAR also has smaller error variance than the others as indicated by its shorter box-plot. The additional `$_+$' sign in the algorithm names (for example, RAAR$_+$) indicates that the algorithms in addition know the amplitude $\mathcal{A}$, i.e., they are applied to the more informative feasibility model (\ref{PR3'}) instead of (\ref{PR3}). The additional information of the amplitude $\mathcal{A}$ improves the performance of every solution method. In this case, we also observe the improved performance of DRAP$_+$ and RAAR$_+$ over VAM$_+$, with average relative RMS errors $4.68$\%, $4.69$\% and $6.82$\%, respectively. The RAAR$_+$ algorithm also has the smallest error variance.}
\label{fig:accuracy comparison}
\end{figure}

We first analyse the performance of the SAM, VAM (equivalently, AP), DR, KM-DR, HPR, RAAR, RRR and DRAP algorithms for solving the feasibility problem (\ref{PR3}), for which recall that the amplitude $\mathcal{A}$ is assumed unknown to the algorithms.
As the DR, KM-DR, HPR and RRR algorithms are clearly outperformed by the RAAR and DRAP methods, we chose to skip their results for brevity.
Table \ref{tbl:beta and iter} shows the number of iterations (the second row), the tuning parameter $\beta$ (the third row) of the algorithms, and the averaged RMS errors over the 75 experiments (the last row).
Due to the extrapolation feature of RAAR and DRAP, each experiment with them is also followed by an averaging process of 20 iterations of alternating projection, indicated by the term `$+\,20$' in the second row of Table \ref{tbl:beta and iter}.
Figure \ref{fig:accuracy comparison} shows the improved performance in terms of accuracy of RAAR and also DRAP over SAM and VAM.
The RAAR algorithm yields phase retrieval with the smallest RMS error on average, $5.98$\% compared to $8.47$\% of SAM, $7.69$\% of VAM and $6.14$\% of DRAP as shown in the last row of Table \ref{tbl:beta and iter}.
RAAR also has smaller error variance than the others as indicated by its shorter box-plot in Figure \ref{fig:accuracy comparison}.
In terms of computational complexity, RAAR and DRAP (50 iterations) are much more efficient than VAM (100 iterations) as shown in Table \ref{tbl:beta and iter} (the second row).
Note that SAM making use of the scalar PSF model has about six times lower complexity per iteration than the other methods; however, this advantage is often dominated by the disadvantage of model deviations for high-NA phase retrieval.

We consider the same 75 high-NA phase retrieval examples as above, but the amplitude $\mathcal{A}$ is now assumed known.
The tighter feasibility model (\ref{PR3'}) then comes into play in place of (\ref{PR3}).
In this section, the algorithms applied to (\ref{PR3'}) will be indicated by the additional `$_+$' sign in their names (for example, RAAR$_+$) to distinguish with them selves for solving (\ref{PR3}).
We analyse the performance of the VAM$_+$ \cite{NguSolVer20} (equivalently, AP$_+$), DR$_+$, KM-DR$_+$, HPR$_+$, RAAR$_+$, RRR$_+$ and DRAP$_+$ algorithms for solving (\ref{PR3'}).
For the same reason as for solving (\ref{PR3}), we chose to skip the phase retrieval results of DR$_+$, KM-DR$_+$, HPR$_+$ and RRR$_+$ for brevity.
The additional information of $\mathcal{A}$ clearly improves the performance of every solution method as shown by Figure \ref{fig:accuracy comparison}, which also demonstrates the improved performance of DRAP$_+$ and RAAR$_+$ over VAM$_+$, with average relative RMS errors $4.68$\%, $4.69$\% and $6.82$\%, respectively.
The RAAR$_+$ algorithm has the smallest error variance.
\bigskip

\noindent\textbf{Funding.} This project has received funding from the ECSEL Joint Undertaking (JU) under grant agreement No. 826589. The JU receives support from the European Union's Horizon 2020 research and innovation programme and Netherlands, Belgium, Germany, France, Italy, Austria, Hungary, Romania, Sweden and Israel.
Russell Luke was supported in part by the Deutsche Forschungsgemeinschaft (DFG, German Research Foundation) – Project-ID 432680300 – SFB 1456 and Project-ID Project-ID LU 1702/1-1.
\bigskip

\noindent
\textbf{Disclosures.} The authors declare no conflicts of interest.

\section*{References}

\bibliographystyle{plain}
\bibliography{shortbib,master_citations}
\end{document}